\documentclass{amsproc}

\usepackage{amsmath,amssymb}
\newtheorem{theorem}{Theorem}[section]
\newtheorem{lemm}[theorem]{Lemma}
\newtheorem{prop}[theorem]{Proposition}

\theoremstyle{definition}
\newtheorem{defi}[theorem]{Definition}
\newtheorem{example}[theorem]{Example}

\newtheorem{coro}[theorem]{Corollary}
\theoremstyle{remark}
\newtheorem{remark}[theorem]{Remark}

\numberwithin{equation}{section}

\def\lg{\langle}
\def\rg{\rangle}

\def\pa{\partial}

\def\al{\alpha}
\def\be{\beta}
\def\ga{\gamma}

\def\ep{\epsilon}

\def\dim{\hbox{dim}}

\def\a{\alpha}

\newfont{\df}{eufm10}

\def\ep{\epsilon}

\def\dim{\hbox{\rm dim}\,}

\begin{document}

\title[Quantizations of Cartan type $S$ Lie algebras]
{Quantizations of generalized Cartan type $S$ Lie algebras \\ and of
the special algebra $\mathbf{S}(n;\underline{1})$ in the modular
case$^\star$}

\author[Hu]{Naihong Hu$^\star$}
\address{Department of Mathematics, East China Normal University,
Min Hang Campus, Dong Chuan Road 500, Shanghai 200241, PR China}
\email{nhhu@math.ecnu.edu.cn}
\thanks{$^\star$N.H.,
supported in part by the NNSF (Grants 10431040, 10728102), the
PCSIRT from the MOE of China, the National/Shanghai Leading Academic
Discipline Project (Project Number: B407).}
\thanks{$^*$X.L., supported in part by the Nankai Research-encouraging Fund for the PhD-Teachers and a fund from LPMC}

\author[Wang]{Xiuling Wang$^*$}
\address{School of Mathematical Sciences and  LPMC, Nankai University,
Tianjin  300071, PR China}\email{xiulingwang@nankai.edu.cn}

\date{Version on Apr. 6nd, 2006}

\begin{abstract}
The generalized Cartan type $\mathbf{S}$ Lie algebras in char $0$
with the Lie bialgebra structures involved are quantized, where the
Drinfel'd twist we used is proved to be a variation of the Jordanian
twist. As the passage from char $0$ to char $p$, their quantization
integral forms are given. By the modular reduction and base changes,
we obtain certain quantizations of the restricted universal
enveloping algebra $\mathbf u(\mathbf{S}(n;\underline{1}))$ (for the
Cartan type simple modular restricted Lie algebra
$\mathbf{S}(n;\underline{1})$ of $\mathbf{S}$ type). They are new
Hopf algebras  of truncated $p$-polynomial noncommutative and
noncocommutative deformation of dimension $p^{1+(n-1)(p^n-1)}$,
which contain the well-known Radford algebra (\cite{DR}) as a Hopf
subalgebra. As a by-product, we also get some Jordanian
quantizations for $\mathfrak {sl}_n$, which are induced from those
horizontal quantizations of $\mathbf S(n;\underline1)$.
\end{abstract}
\maketitle

In Hopf algebra or quantum group theory, two standard methods to
construct new bialgebras from old ones are by twisting the product
by a $2$-cocycle but keeping the coproduct unchanged, and by
twisting the coproduct by a Drinfel'd twist but preserving the
product. Constructing quantizations of Lie bialgebras is an
important approach to producing new quantum groups (see \cite{D,ES}
and references therein). In papers \cite{EK1, EK2}, Etingof-Kazhdan
showed the existence of a universal quantization for Lie bialgebras
by constructing  a quantization functor. Enriquez-Halbout showed
that any coboundary  Lie bialgebra, in principle, can be quantized
via a certain Etingof-Kazhdan quantization functor (see \cite{EH}).
The Lie bialgebras they considered (including finite- and
infinite-dimensional ones) are those classes of the Lie algebras
defined by generalized Cartan matrices. However, for another
important class of the Cartan type Lie algebras defined by
differential operators, they are lack of adequate attention in the
literature. In 2004, Grunspan \cite{CG} obtained the quantization of
the (infinite-dimensional) Witt algebra $W$ in characteristic $0$
using the twist found by Giaquinto-Zhang (\cite{AJ}), but his way
didn't work for the quantum version (\cite{CG}) of its simple
modular Witt algebra $\mathbf{W}(1;\underline{1})$ in characteristic
$p$. In 2005, Song-Su (\cite{GY}) determined some coboundary
triangular Lie bialgebra structures (for the definition, see p. 28,
\cite{ES}) on the Lie algebras of the generalized Witt type. The
authors (\cite{HW}) obtained the quantizations both for the
generalized-Witt algebra $\mathbf{W}$ in characteristic $0$ and for
the Jacobson-Witt algebra $\mathbf{W}(n;\underline{1})$ in
characteristic $p$, where in the rank $1$ case, we recovered the
Grunspan's work in characteristic $0$, and gave the required quantum
version in characteristic $p$.

\smallskip
In the present paper, we continue to treat with the same questions
both for the generalized Cartan type $\mathbf{S}$ Lie algebras in
characteristic $0$ (for the definition, see \cite{DZ1}) and for the
restricted simple special algebras $\mathbf{S}(n;\underline{1})$ in
the modular case (for the definition, see \cite{H}, \cite{HR}).

\smallskip
As we known, the construction of Drinfeld twists is difficult. Only
a few of twists in explicit forms have been known for a long time
(see \cite{AJ,KL,VO,NR} etc.). In this paper, we start an explicit
Drinfel'd twist due to \cite{AJ} and \cite{CG}, which, we found
recently, is essentially a variation of the Jordanian twist used by
Kulish et al (see \cite{KL}, etc.). For this fact, we provide a
strict proof in Remark 1.10. Similar to \cite{CG} and \cite{HW}, we
quantize the triangular Lie bialgebra structures on the generalized
Cartan type $\mathbf{S}$ Lie algebra in characteristic $0$ (The
existence of triangular Lie bialgebra structures on it is due to a
result of \cite{M}). The process depends on the construction of
Drinfel'd twists which, up to integral scalars, are controlled by
the classical Yang-Baxter $r$-matrix. To study the modular case,
what we discuss first is about the arithmetic property of
quantizations to work out their quantization integral forms. To this
end, we have to work over the so-called {\it ``positive"} part
subalgebra $\mathbf{S}^+$ of the generalized Cartan type
$\mathbf{S}$ Lie (shifted) algebra $x^{\eta}\mathbf{S}$ (where
$\eta=-\underline{1}$). This is one of the crucial technical points
here. It is an infinite-dimensional simple Lie algebra when defined
over a field of characteristic $0$, while, defined over a field of
characteristic $p\,$, it contains a maximal ideal
$J_{\underline{1}}$ and the corresponding quotient is exactly the
algebra $\mathbf{S}'(n;\underline{1})$. Its derived subalgebra
$\mathbf{S}(n;\underline{1})=\mathbf{S}'(n;\underline{1})^{(1)}$ is
a Cartan type restricted simple modular Lie algebra of $\mathbf{S}$
type. Secondly, in order to yield the {\it expected}
finite-dimensional quantizations of the restricted universal
enveloping algebra of the special algebra
$\mathbf{S}(n;\underline{1})$, we need to carry out the modular
reduction process: {\it modulo $p$ reduction} and {\it modulo
``$p$-restrictedness" reduction}, among which, we have to take the
suitable {\it base changes}. These are the other two crucial
technical points. Thirdly, in the process of giving the quantization
integral forms for the $\mathbb{Z}$-form $\mathbf{S}_{\mathbb Z}^+$
in characteristic $0$, we find that there exist $n(n-1)$ the
so-called {\it basic Drinfel'd twists}, which can afford many more
Drinfel'd twists (Corollary 5.2). Furthermore, we investigate the
twisted structures arisen by these twists. Note that these Hopf
algebras contain the well-known Radford algebra (\cite{DR}) as a
Hopf subalgebra. Our work gets a new class of noncommutative and
noncocommutative finite-dimensional Hopf algebras in characteristic
$p$ (see \cite{T}).

\smallskip
The article is organized as follows: In Section 1, we collect some
definitions and lemmas which are useful for later use. In Section 2,
we quantize explicitly Lie bialgebra structures of generalized
Cartan type $\mathbf{S}$ Lie algebra $x^{\eta}\mathbf{S}$ by the
\textit{basic Drinfeld twists} ({\it in vertical}), and obtain
$n(n-1)$ new quantization integral forms for $\mathbf{S}_{\mathbb
Z}^+$ in characteristic $0$. We use this fact to equip the
restricted universal enveloping algebra of the special algebra
$\mathbf{S}(n;\underline{1})$ with noncommutative and
noncocommutative Hopf algebra structures by the modular reduction
and the base changes in Section 3. In Section 4, considering some
products of pointwise different basic Drinfel'd twists, we can get
new quantization integral forms for $\mathbf{S}_{\mathbb Z}^+$ in
characteristic $0$, which, via modulo $p$ reduction and
``$p$-restrictedness" reduction, together with two steps of base
changes, eventually leads to new Hopf algebras of dimension
$p^{1{+}(n-1)(p^n-1)}$ with indeterminate $t$ or of dimension
$p^{(n-1)(p^n-1)}$ with specializing $t$ into a scalar in
$\mathcal{K}$ in characteristic $p$. In Section 5, using the {\it
horizontal} twists (different from the {\it vertical} ones in
Section 3), we get some new quantizations of horizontal type of
$\mathbf u(\mathbf{S}(n;\underline{1}))$, which contain some
quantizations of the Lie algebra $\mathfrak{sl}_n$ derived by the
Jordanian twists (cf. \cite{KLS}).

\section{Preliminaries}
\medskip
\subsection{Generalized Cartan type S Lie  algebra and its Lie bialgebra structure}
Let $\mathbb{F}$ be a field with $\text{char}(\mathbb{F})=0$ and
$n>0$. Let $\mathbb{Q}_n=\mathbb{F}[x^{\pm1}_1,\cdots,x^{\pm1}_n]$
be a Laurent polynomial algebra and $\partial_i$ coincides with the
degree operator $x_i\frac{\partial}{\partial x_i}$.
 Set
$T=\bigoplus_{i=1}^n{\mathbb{Z}}\partial_i$, and
$x^{\alpha}=x_1^{\alpha_1} \cdots x_n^{\alpha_n}$ for
$\alpha=(\alpha_1,\cdots,\alpha_n)
 \in \mathbb{Z}^{n}$.

Denote $\mathbf{W}=\mathbb{Q}_n\otimes_{\mathbb{Z}}
T=\text{Span}_{\mathbb{F}}\{\,x^{\alpha}\partial \mid \alpha \in
\mathbb{Z}^{n},\partial\in T\}$, where we set $x^{\alpha}\partial\
=x^{\alpha}\otimes
\partial$ for short. Then $\mathbf{W}=\text{Der}_{\mathbb{F}}(\mathbb{Q}_n)$
is a Lie algebra of
generalized-Witt type (see \cite{DZ}) under the following bracket
$$[\,x^{\alpha}\partial,\,x^{\beta}\partial'\,]
 =x^{\alpha+\beta}\bigl(\,\partial(\beta)\partial'-\partial'(\alpha)\partial\,\bigr),
 \qquad\forall\
 \alpha,\, \beta\in \mathbb{Z}^{n};\ \partial,\, \partial' \in T,$$
where
$\partial(\beta)=\lg\partial,\beta\rg=\lg\beta,\partial\rg=\sum\limits_{i=1}^{n}a_i\beta_i\in\mathbb{Z}$
for $\partial=\sum\limits_{i=1}^{n}a_i\partial_i \in T$ and
$\beta=(\beta_1,\cdots,\beta_n) \in \mathbb{Z}^{n}$. The bilinear
map $\lg\cdot,\cdot\rg: \,T\times {\mathbb{Z}}^n\longrightarrow
{\mathbb{Z}}$ is non-degenerate in the sense that
\begin{gather*}
\partial(\al)=\lg\partial,\al\rg=0 \quad (\forall\;\partial \in
T),\ \Longrightarrow \al=0,\\
\partial(\al)=\lg\partial,\al\rg=0 \quad  (\forall\;\al \in
\mathbb{Z}^n),\ \Longrightarrow
\partial=0.
\end{gather*}
$\mathbf{W}$ is an infinite dimensional simple Lie algebra over
$\mathbb{F}$  (see \cite{DZ}).

We recall that the $divergence$ (cf. \cite{DZ1}) $\mathrm{div}:
\mathbf{W}\longrightarrow \mathbb{Q}_n$ is the $\mathbb{F}$-linear
map such that
$$
\mathrm{div}(x^{\al}\partial)=\partial(x^{\al})=\partial(\al)x^{\al},\quad
\textit{for } \al \in \mathbb{Z}^n,\ \partial \in T.\leqno (1)
$$
The $divergence$ has the following two properties:
\begin{gather*}
\mathrm{div} ([u,v])=u\cdot \mathrm{div}(v)-v\cdot \mathrm{div}(u),\tag{2}\\
\mathrm{div}(fw)=f\,\mathrm{div}(w)+w\cdot f, \tag{3}
\end{gather*}
for $u,v,w \in \mathbf{W}, f \in \mathbb{Q}_n$. In view of (2), the
subspace
$$\widetilde{\mathbf{S}}=\mathrm{Ker}(\mathrm{div})$$ is a subalgebra of $\mathbf{W}$.

The Lie algebra $\mathbf{W}$ is $\mathbb{Z}^n$-graded, whose
homogeneous components are
$$
\mathbf{W}_{\al}:=x^{\al}T, \qquad \al
\in\mathbb{Z}^n.
$$
The \textit{divergence} $\mathrm{div}:
\mathbf{W}\longrightarrow \mathbb{Q}_n$ is a derivation of degree
$0$. Hence, its kernel is a homogeneous subalgebra of $\mathbf{W}$.
So we have
$$\widetilde{\mathbf{S}}=\bigoplus_{\al \in \mathbb{Z}^n}
\widetilde{\mathbf{S}}_{\al}, \quad
\widetilde{\mathbf{S}}_{\al}:=\widetilde{\mathbf{S}}\cap
\mathbf{W}_{\al}.$$ For each $\al \in \mathbb{Z}^n,$ let
$\hat{\al}: T\rightarrow \mathbb{F}$ be the corresponding linear
function defined by
$\hat{\al}(\partial)=\lg\partial,\al\rg=\partial(\al)$. We have
$$
\widetilde{\mathbf{S}}_{\al}=x^{\al}T_{\al}, \quad\textit{and } \
T_{\al}=\mathrm{Ker}(\hat{\al}).
$$
The algebra $\widetilde{\mathbf{S}}$ is not simple, but its derived
subalgebra $\mathbf{S}=(\widetilde{\mathbf{S}})'$ is simple,
assuming only that $\mathrm{dim}  T\geq 2$. According to Proposition
3.1 \cite{DZ1}, we have $\mathbf{S}=\bigoplus_{\al\neq
0}\widetilde{\mathbf{S}}_{\al}$. More generally, it turns out that
the shifted spaces $x^{\eta}\mathbf{S},\, \eta \in
\mathbb{Z}^n-\{0\}$, are simple subalgebras of $\mathbf{W}$ if
$\mathrm{dim}  T\geq 3$. We refer to the simple Lie algebras
$x^{\eta}\mathbf{S}$ as the Lie algebras of {\it generalized  Cartan
type} $\mathbf{S}$ (see \cite{DZ1}). The Lie algebra
$x^{\eta}\mathbf{S}$ is $\mathbb{Z}^n$-graded with
$x^{\al}T_{\al-\eta} \ (\al\neq \eta)$, as its homogeneous component
of degree $\al$, while its homogeneous component of degree $\eta$ is
$0$.

Throughout this paper, we assume that $\eta\neq 0$,
$\eta_k=\eta_{k'}$.

Take two distinguished elements $h=\pa_k-\pa_{k'}, e=x^{\gamma}\pa_0
\in  x^{\eta}\mathbf{S}$ such that  $[h,e]=e$ where $1\leq k\neq k'
\leq n$. It is easy to see that $\pa_0(\gamma-\eta)=0$, and
$\gamma_k-\gamma_{k'}=1$. Using a result of  \cite{M}, we have the
following
\begin{prop}
There is a triangular Lie bialgebra structure on
$x^{\eta}\mathbf{S}$ $(\eta\neq 0, \ \eta_k=\eta_{k'})$  given by
the classical Yang-Baxter $r$-matrix
$$
r:=(\pa_k-\pa_{k'})\otimes x^{\gamma}\partial_0-x^{\gamma}\partial_0
\otimes (\partial_k-\pa_{k'}),\quad \forall \;\pa_{k'},\,\pa_k \in
T, \ \gamma \in \mathbb{Z}^{n},
$$
where $\gamma_k-\gamma_{k'}=1$, $\pa_0(\gamma)=\pa_0(\eta)$ and
$[\,\partial_k-\pa_{k'},
x^\gamma\partial_0\,]=x^\gamma\partial_0$.\hfill\qed
\end{prop}

\subsection{Generalized Cartan type $\mathbf{S}$ Lie subalgebra $\mathbf{S}^+$}
Denote $D_i=\frac{\partial}{\partial x_i}$. Set
$\mathbf{W}^+:=\text{Span}_{\mathcal{K}}\{x^{\alpha} D_i\mid
\alpha\in\mathbb{Z}_+^n, 1\le i\le n\}$, where $\mathbb{Z}_+$ is the
set of non-negative integers. Then
$\mathbf{W}^+=\text{Der}_{\mathcal{K}}(\mathcal{K}[x_1,\cdots,x_n])$
is the derivation Lie algebra of polynomial ring
$\mathcal{K}[x_1,\cdots,x_n]$, which, via the identification $x^\al
D_i$ with $x^{\alpha-\epsilon_i}
\partial_i$ (here $\alpha-\epsilon_i\in\mathbb{Z}^n$, $\epsilon_i=(\delta_{1i},\cdots,\delta_{ni})$),
can be viewed as a Lie subalgebra (the ``positive" part) of the
generalized-Witt algebra $\mathbf{W}$ over a field $\mathcal{K}$.

For $X=\sum_{i=1}^{n}a_iD_i \in \mathbf{W}$, we define
$\mathrm{Div}(X)=\sum_{i=1}^{n}D_i(a_i)$ as usual. Note that
$\mathrm{div}(X)=\sum_{i=1}^{n}x_iD_i(x_i^{-1}a_i)$ (since
$\pa_i=x_iD_i$). Thus we have $\mathrm{div}(x_1\cdots x_n
X)=x_1\cdots x_n \mathrm{Div}(X)$. This means that
$X\in\mathrm{Ker}(\mathrm{Div})$ if and only if
$x^{\underline{1}}X\in \widetilde {\mathbf{S}}$, and if and only if
$X\in x^{-\underline{1}}\widetilde{\mathbf{S}}$, where
$\underline{1}=\epsilon_1+\cdots+\epsilon_n$.

Set $\mathbf{S}^+:=\mathrm{Ker}(\mathrm{Div})\bigcap\mathbf{W}^+$,
then we have
$\mathbf{S}^+=(x^{-\underline{1}}\mathbf{S})\bigcap\mathbf{W}^+$
since $\mathbf{S}=\bigoplus_{\al\neq 0}\widetilde{\mathbf{S}}_{\al}$
and
$x^{-\underline{1}}\widetilde{\mathbf{S}}_0\bigcap\mathbf{W}^+=0$
(where $\widetilde{\mathbf{S}}_0=T$), which is a subalgebra of
$\mathbf{W}^+$. Note that
$\{\,\alpha_nx^{\alpha-\epsilon_n}D_i-\alpha_ix^{\alpha-\epsilon_i}D_n\mid
\al\in\mathbb{Z}^n_+, 1\le i< n\,\}$ is a basis of $\mathbf{S}^+$,
where
$\alpha_nx^{\alpha-\epsilon_n}D_i-\alpha_ix^{\alpha-\epsilon_i}D_n
=x^{\al-\ep_i-\ep_n}(\al_n\pa_i-\al_i\pa_n)\in
x^{\al-\ep_i-\ep_n}T_{\al-\ep_i-\ep_n+\underline{1}}$ indicates once
again that $\mathbf{S}^+$ is indeed a subalgebra of
$x^{-\underline{1}}\mathbf{S}$ since $\pa_i=x_iD_i$.


\subsection{The special algebra $\mathbf{S}(n;\underline{1})$}

Assume now that $\text{char}(\mathcal{K})=p$, then by definition,
the Jacobson-Witt algebra $\mathbf{W}(n;\underline{1})$ is a
restricted simple Lie algebra over a field $\mathcal{K}$. Its
structure of $p$-Lie algebra is given by $D^{[p]}=D^p,\; \forall\, D
\in \mathbf{W}(n;\underline{1})$ with a basis $\{\,x^{(\alpha)}D_j
\mid 1\leq j\leq n, \ 0 \leq \alpha \leq \tau \}$, where
$\tau=(p{-}1,\cdots,p{-}1) \in \mathbb{N}^n$;
$\epsilon_i=(\delta_{1i},\cdots,\delta_{ni})$ such that
$x^{(\epsilon_i)}=x_i$ and $D_j(x_i)=\delta_{ij}$; and $\mathcal
O(n;\underline{1}):=\{\,x^{(\al)}\mid 0 \leq \alpha \leq \tau \}$ is
the restricted divided power algebra with
$x^{(\al)}x^{(\be)}=\binom{\al{+}\be}{\al}\,x^{(\al{+}\be)}$ and a
convention: $x^{(\alpha)}=0$ if $\alpha$ has a component
$\alpha_j<0$ or $\geq p$, where
$\binom{\al{+}\be}{\al}:=\prod_{i=1}^n\binom{\al_i{+}\be_i}{\al_i}$.
Note that $\mathcal O(n;\underline{1})$ is $\mathbb{Z}$-graded by
$\mathcal
O(n;\underline{1})_i:=\textrm{Span}_{\mathcal{K}}\{\,x^{(\al)}\mid 0
\leq \alpha \leq \tau, |\a|=i \}$, where $|\a|=\sum_{j=1}^{n}\a_j$.
Moreover, $\mathbf{W}(n;\underline{1})$ is isomorphic to
$\text{Der}_{\mathcal{K}}(\mathcal O(n;\underline{1}))$ and inherits
a gradation from $\mathcal O(n;\underline{1})$ by means of
$\mathbf{W}(n;\underline{1})_i=\sum_{j=1}^{n}\mathcal
O(n;\underline{1})_{i+1}D_j$. Then the subspace
$$\mathbf{S}'(n;\underline{1})=\{E \in \mathbf{W}(n;\underline{1}) \mid \mathrm{Div}(E)=0\}$$
is a subalgebra of  $\mathbf{W}(n;\underline{1})$.

Its derived subalgebra
$\mathbf{S}(n;\underline{1})=\mathbf{S}'(n;\underline{1})^{(1)}$ is
called \textit{the special algebra}. Then
$\mathbf{S}(n;\underline{1})=\bigoplus_{i=-1}^{s}\mathbf{S}(n;\underline{1})
\cap\mathbf{W}(n;\underline{1})_i$ is graded with $s=|\tau|-2$.
Recall the mappings $D_{ij}: \mathcal
O(n;\underline{1})\longrightarrow \mathbf{W}(n;\underline{1})$,
$D_{ij}(f)=D_j(f)D_i-D_i(f)D_j$ for $f \in \mathcal
O(n;\underline{1})$. Note that $D_{ii}=0$ and $D_{ij}=-D_{ji},\
1\leq i,\,j \leq n$. Then by Lemma 4.2.2 \cite{H},
$$\mathbf{S}(n;\underline{1})=\text{Span}_{\mathcal{K}}\{D_{in}(f) \mid f \in
\mathcal O(n;\underline{1}), 1\leq i< n\}$$ is a $p$-subalgebra of
$\mathbf{W}(n;\underline{1})$ with restricted gradation. Evidently,
we have the following result (see the proof of Theorem 3.7, p.159 in
\cite{HR})
\begin{lemm}
$\mathbf{S}'(n;\underline{1})=\mathbf{S}(n;\underline{1})+\sum\limits_{j=1}^{n}\mathcal{K}
x^{(\tau-(p-1)\epsilon_j)}D_j$. And
$\mathrm{dim}_{\mathcal{K}}\mathbf{S}'(n;\underline{1})$
$=(n-1)p^n+1$,
$\mathrm{dim}_{\mathcal{K}}\mathbf{S}(n;\underline{1})=(n-1)(p^n-1)$.\hfill\qed
\end{lemm}

By definition (cf. \cite{H}), the restricted universal enveloping
algebra $\mathbf{u}(\mathbf{S}(n;\underline{1}))$ is isomorphic to
$U(\mathbf{S}(n;\underline{1}))/I$ where $I$ is the Hopf ideal of
$U(\mathbf{S}(n;\underline{1}))$ generated by
$(D_{ij}(x^{(\epsilon_i+\epsilon_j)}))^p-D_{ij}(x^{(\epsilon_i+\epsilon_j)}),\,
(D_{ij}(x^{(\a)}))^p$ with $\a \neq \epsilon_i+\epsilon_j$ for
$1\leq i<j\leq n$. Since
$\mathrm{dim}_{\mathcal{K}}\mathbf{S}(n;\underline{1})=(n-1)(p^n-1)$,
we have
$\mathrm{dim}_{\mathcal{K}}\mathbf{u}(\mathbf{S}(n;\underline{1}))$
$=p^{(n-1)(p^n-1)}$.

\subsection{A crucial Lemma}
For any element $x$ of a unital $R$-algebra ($R$ a ring) and
     $a \in R$, we set (see \cite{AJ})
$$x_a^{\lg n\rg}:=(x+a)(x+a+1)\cdots(x+a+n-1), \leqno(4)$$
then $x^{\lg n\rg}:=x_0^{\lg n\rg}=\sum_{k=0}^nc(n,k)x^k$ where
$c(n,k)$ is the number of $\pi\in \frak S_n$ with exactly $k$ cycles
(cf. \cite{R}). Given a $\pi\in \frak S_n$, let $c_i=c_i(\pi)$ be
the number of cycles of $\pi$ of length $i$. Note that $n=\sum
ic_i$. Define the type of $\pi$, denoted type $\pi$, to be the
$n$-tuple $\underline{c}=(c_1,\cdots,c_n)$. The total number of
cycles of $\pi$ is denoted $c(\pi)$, so
$c(\pi)=|\,\underline{c}\,|=c_1+\cdots+c_n$. Denote by $\frak
S_n(\underline{c})$ the set of all $\sigma\in \frak S_n$ of type
$\underline{c}$, then $|\frak
S_n(\underline{c})|=n!/1^{c_1}c_1!2^{c_2}c_2!\cdots n^{c_n}c_n!$
(see Proposition 1.3.2 \cite{R}).

We also set
$$x_a^{[n]}:=(x+a)(x+a-1)\cdots(x+a-n+1), \leqno(5)$$
then $x^{[n]}:=x_0^{[n]}=\sum_{k=0}^n s(n,k)x^k$ where
$s(n,k)=(-1)^{n-k}c(n,k)$ is the Stirling number of the first kind.

\begin{lemm} $($\cite{AJ,CG}$)$
For any element $x$ of a unital $\mathbb{F}$-algebra with
$\text{char}(\mathbb{F})=0$, $a, \,b \in \mathbb{F}$ and $r,\, s,\,
t \in \mathbb{Z}$, one has
\begin{gather*}
x_a^{\lg s+t\rg}=x_a^{\lg s\rg}\,x_{a+s}^{\lg t\rg},\tag{6} \\
x_a^{[s+t]}=x_a^{[s]}\,x_{a-s}^{[t]},\tag{7} \\
x_a^{[s]}=x_{a-s+1}^{\lg s\rg},\tag{8} \\
\sum\limits_{s+t=r}\frac{(-1)^t}{s!\,t!}\,x_a^{[s]}\,x_b^{\lg
t\rg}=\dbinom{a{-}b} {r}=\frac{(a{-}b)\cdots(a{-}b{-}r{+}1)}{r!}, \tag{9} \\
\sum\limits_{s+t=r}\frac{(-1)^t}{s!\,t!}\,x_a^{[s]}\,x_{b-s}^{[t]}=\dbinom{a{-}b{+}r{-}1}
{r}=\frac{(a{-}b)\cdots(a{-}b{+}r{-}1)}{r!}.\tag{10}
\end{gather*}
\end{lemm}

\subsection{Quantization by Drinfel'd twists} The following
result is well-known (see \cite{CP,D,ES,NR}, etc.).
\begin{lemm}
Let $(A,m,\iota,\Delta_0,\varepsilon_0,S_0)$ be a Hopf algebra over
a commutative ring. A Drinfel'd twist $\mathcal{F}$ on $A$ is an
invertible element of $A\otimes A$ such that
\begin{gather*}(\mathcal{F}\otimes
1)(\Delta_0\otimes \text{\rm Id})(\mathcal{F})=(1\otimes
\mathcal{F})(\text{\rm Id}\otimes\Delta_0)(\mathcal{F}), \\
(\varepsilon_0\otimes \text{\rm Id})(\mathcal{F})=1=(\text{\rm
Id}\otimes \varepsilon_0)(\mathcal{F}).
\end{gather*} Then,
$w=m(\text{\rm Id}\otimes S_0)(\mathcal{F})$ is invertible in $A$
with $w^{-1}=m(S_0\otimes \text{\rm Id})(\mathcal{F}^{-1})$.

Moreover, if we define $\Delta: \,A\longrightarrow A\otimes A$ and
$S: \,A\longrightarrow A$ by
$$\Delta(a)=\mathcal{F}\Delta_0(a)\mathcal{F}^{-1},
\qquad S=w\,S_0(a)\,w^{-1},$$ then $(A, m, \iota,
\Delta,\varepsilon,S)$ is a new Hopf algebra, called the twisting
of $A$ by the Drinfel'd twist $\mathcal{F}$.
\end{lemm}

Let $\mathbb{F}[[t]]$ be a ring of formal power series over a field
$\mathbb{F}$ with $\text{char}(\mathbb{F})=0$. Assume that $L$ is a
triangular Lie bialgebra over $\mathbb{F}$ with a classical
Yang-Baxter $r$-matrix $r$ (see \cite{D,ES}). Let $U(L)$ denote the
universal enveloping algebra of $L$, with the standard Hopf algebra
structure $(U(L),m,\iota,\Delta_0,\varepsilon_0,S_0)$.

\smallskip
Let us consider {\it the topologically free
$\mathbb{F}[[t]]$-algebra} $U(L)[[t]]$ (for the definition, see p.
4, \cite{ES}), which can be viewed as an associative
$\mathbb{F}$-algebra of formal power series with coefficients in
$U(L)$. Naturally, $U(L)[[t]]$ equips with an induced Hopf algebra
structure arising from that on $U(L)$ (via the coefficient ring
extension), by abuse of notation, denoted still by
$(U(L)[[t]],m,\iota,\Delta_0,\varepsilon_0,S_0)$.

\begin{defi} (\cite{HW})
For a triangular Lie bialgebra $L$ over $\mathbb{F}$ with
$\text{char}(\mathbb{F})=0$, $U(L)[[t]]$ is called {\it a
quantization of $U(L)$ by a Drinfel'd twist} $\mathcal{F}$ over
$U(L)[[t]]$ if $U(L)[[t]]/tU(L)[[t]]\cong U(L)$, and $\mathcal{F}$
is determined by its $r$-matrix $r$ (namely, its Lie bialgebra
structure).
\end{defi}

\subsection{Construction of Drinfel'd twists}
Let $L$ be a Lie algebra containing linearly independent elements
$h$ and $e$ satisfying $[h,\,e]=e$, then the classical Yang-Baxter
$r$-matrix $r=h\otimes e-e\otimes h$ equips $L$ with the structure
of triangular coboundary Lie bialgebra (see \cite{M}). To describe a
quantization of $U(L)$ by a Drinfel'd twist $\mathcal{F}$ over
$U(L)[[t]]$, we need an explicit construction for such a Drinfel'd
twist. In what follows, we shall see that such a twist depends upon
the choice of {\it two distinguished elements} $h,\,e$ arising from
its $r$-matrix $r$.

Recall the following results proved in \cite{CG} and \cite{HW}. Note
that $h$ and $e$ satisfy the following equalities:
$$e^s\cdot h_a^{[m]}=h_{a-s}^{[m]}\cdot e^s, \leqno{(11)}$$
$$e^s\cdot h_a^{\lg m\rg}=h_{a-s}^{\lg m\rg}\cdot e^s, \leqno{(12)}$$
where $m,\,s$ are non-negative integers, $a \in \mathbb{F}$.

For $a \in \mathbb{F}$, following \cite{CG}, we set
\begin{gather*}
\mathcal{F}_a=\sum\limits_{r=0}^{\infty}\frac{(-1)^r}{r!}h_a^{[r]}\otimes
e^rt^r,\qquad F_a=\sum\limits_{r=0}^{\infty}\frac{1}{r!}h_a^{\lg
r\rg}\otimes
e^rt^r,\\
u_a=m\cdot(S_0\otimes \text{\rm Id})(F_a),\qquad\quad
v_a=m\cdot(\text{\rm Id}\otimes S_0)(\mathcal{F}_a).
\end{gather*}

Write $\mathcal{F}=\mathcal{F}_0,\, F=F_0,\,u=u_0,\,v=v_0$.

Since $S_0(h_a^{\lg r\rg})=(-1)^rh_{-a}^{[r]}$ and
$S_0(e^r)=(-1)^re^r$, one has
$$
v_a=\sum\limits_{r=0}^{\infty}\frac{1}{r!}h_a^{[r]}
e^rt^r, \quad
u_b=\sum\limits_{r=0}^{\infty}\frac{(-1)^r}{r!}h_{-b}^{[r]}
e^rt^r.
$$
\begin{lemm} $($\cite{CG}$)$
For $a,\, b \in \mathbb{F}$, one has
$$
\mathcal{F}_a F_b=1\otimes(1-et)^{a-b} \quad\text{and }\quad v_a
u_b=(1-et)^{-(a+b)}.
$$
\end{lemm}

\begin{coro} $($\cite{CG}$)$
For $a \in \mathbb{F}$, $\mathcal{F}_a$ and $u_a$ are invertible
with $\mathcal{F}_a^{-1}=F_a$ and $u_a^{-1}=v_{-a}$.  In particular,
$\mathcal{F}^{-1}=F$ and $u^{-1}=v$.
\end{coro}
\begin{lemm} $($\cite{CG}$)$ For any positive integers $r$, we have
$$\Delta_0(h^{[r]})=\sum\limits_{i=0}^r \dbinom{r}{i}h^{[i]}\otimes
h^{[r-i]}.$$  Furthermore, $\Delta_0(h^{[r]})=\sum\limits_{i=0}^r
\dbinom{r}{i}h^{[i]}_{-s}\otimes h^{[r-i]}_s$ for any $s \in
\mathbb{F}$.
\end{lemm}

\begin{prop} $($\cite{CG}$)$ If a Lie algebra $L$ contains a $2$-dimensional solvable
Lie subalgebra with a basis $\{h,\,e\}$ satisfying $[h,\,e]=e$, then
$\mathcal{F}=\sum\limits_{r=0}^{\infty}\frac{(-1)^r}{r!}h^{[r]}\otimes
e^rt^r$ is a Drinfel'd twist on $U(L)[[t]]$.
\end{prop}

\begin{remark} Recently, we observed that Kulish et al earlier used the so-called
{\it Jordanian twist} (see \cite{KL}) with the two-dimensional
carrier subalgebra $B(2)$ such that $[H,E]=E$, defined by the
canonical twisting element
$$\mathcal{F}_{\mathcal{J}}^{c}=\mathrm{exp}(H\otimes\sigma(t)), \quad
\sigma(t)=\mathrm{ln}(1+Et),$$ where
$\exp(X)=\sum_{i=0}^{\infty}\frac{X^n}{n!}$ and
$\ln(1+X)=\sum_{n=1}^{\infty}\frac{(-1)^{n+1}}{n}X^n$.

Expanding it, we get
\begin{equation*}
\begin{split}
\exp(H\otimes\sigma(t))
&=\exp\Bigl(\sum_{n=1}^{\infty}\frac{(-1)^{n+1}}{n}H\otimes (Et)^n\Bigr)\\
&=\prod_{n\ge 1}\Bigl(\sum_{\ell\ge
0}\frac{(-1)^{(n{+}1)\ell}}{n^\ell \ell!}H^\ell\otimes
(Et)^{n\ell}\Bigr)
\\
&=\sum_{n\ge1}\sum_{c_1,\cdots,c_n\ge0}\frac{(-1)^{c_1+2c_2+\cdots+nc_n-|\,\underline{c}\,|}}
{c_1!\cdots c_n!1^{c_1}2^{c_2}\cdots
n^{c_n}}H^{|\,\underline{c}\,|}\otimes
(Et)^{c_1+2c_2+\cdots+nc_n} \\
&=\sum_{n\ge0}\Bigl(\sum_{\underline{c}}\frac{(-1)^{n-|\,\underline{c}\,|}|\,\frak
S_n(\underline{c})\,|} {n!}H^{|\,\underline{c}\,|}\Bigr)\otimes
(Et)^n        \\
&=\sum_{n\ge0}\left(\sum_{k=0}^n\frac{(-1)^{n-k}c(n,k)}
{n!}H^k\right)\otimes
(Et)^n        \\
&=\sum_{n=0}^{\infty}\frac{1}{n!}H^{[n]}\otimes E^nt^n, \qquad
\end{split}
\end{equation*}
where we set $n=c_1+2c_2+\cdots+nc_n$,
$c(n,k)=\sum_{|\,\underline{c}\,|=k}|\,\frak S_n(\underline{c})\,|$.
So
$$(\mathcal{F}_{\mathcal{J}}^{c})^{-1}=\exp((-H)\otimes\sigma(t))=\sum_{r=0}^{\infty}
\frac{1}{r!}(-H)^{[r]}\otimes
E^rt^r=\sum_{r=0}^{\infty}\frac{(-1)^r}{r!}H^{\langle
r\rangle}\otimes E^rt^r.$$

Consequently, we can rewrite the twist $\mathcal{F}$ in
Proposition 1.9 as
$$\mathcal{F}=\sum_{r=0}^{\infty}
\frac{(-1)^r}{r!}H^{[r]}\otimes E^rt^r=\exp(H\otimes\sigma'(t)),
\quad \sigma'(t)=\ln(1-Et),$$ where $[H,-E]=-E$. So there is no
difference between the twists $\mathcal{F}$ and
$\mathcal{F}_{\mathcal{J}}^{c}$. They are essentially the same up to
an isomorphism on the carrier subalgebra $B(2)$.
\end{remark}

\section{Quantization of Lie bialgebra of generalized Cartan type $\mathbf{S}$}
In this section, we explicitly quantize the Lie bialgebras
$x^{\eta}\mathbf{S}$ of generalized Cartan type $\mathbf{S}$ by the
twist given in Proposition 1.9.
\subsection{Some commutative relations in $U(x^{\eta}\mathbf{S})$}
For the universal enveloping algebra $U(x^{\eta}\mathbf{S})$ of the
generalized Cartan type $\mathbf{S}$ Lie algebra
$x^{\eta}\mathbf{S}$ over $\mathbb{F}$, we need to do some necessary
calculations, which are important to the quantizations  of Lie
bialgebra structure of $x^{\eta}\mathbf{S}$ in the sequel.

\begin{lemm} Fix two distinguished
elements $h:=\partial_k{-}\pa_{k'}$, $e:=x^{\gamma}\partial_0 \in
x^{\ga} T_{\ga-\eta}$ with $\gamma_k-\gamma_{k'}=1$ for
$x^{\eta}\mathbf{S}$. For $a \in \mathbb{F}$, $x^{\alpha}\partial
\in x^{\al}T_{\al-\eta},\, x^{\beta}\partial' \in
x^{\beta}T_{\beta-\eta}$,  $m$ is non-negative integer, the
following equalities hold in $U(x^{\eta}\mathbf{S}):$
\begin{gather*}
x^{\alpha}\partial\cdot h_a^{[m]}=h_{a+(\a_{k'}-\al_k)}^{[m]}\cdot
x^{\alpha}\partial,
\tag{13}\\
x^{\alpha}\partial\cdot h_a^{\lg m\rg}=h_{a+(\a_{k'}-\al_k)}^{\lg
m\rg}\cdot x^{\alpha}\partial, \tag{14}\\
x^{\alpha}\partial\cdot(x^{\beta}\partial')^m=\sum\limits_{\ell=0}^m({-}1)^\ell\dbinom{m}
{\ell} (x^{\beta}\partial')^{m{-}\ell}\cdot
x^{\alpha{+}\ell\beta}\Bigl(a_\ell
\partial-b_\ell\partial'\Bigr), \tag{15}
\end{gather*}
where $a_\ell=\prod\limits_{j=0}^{\ell-1}\partial'(\alpha{+}j\beta)=
\prod\limits_{j=0}^{\ell-1}\partial'(\alpha{+}j\eta)$,
$b_\ell=\ell\,\partial(\beta)a_{\ell{-}1}$, and set $a_0=1$,
$b_0=0$.
\end{lemm}
\begin{proof}
One has (13) and (14) by using induction on $m$.

Formula (15) is a consequence of the fact (see Proposition 1.3 (4),
\cite{HR}) that for any elements $a,\,c$ in an associative algebra,
one has
$$c\,a^m=\sum_{\ell=0}^m(-1)^{\ell}\dbinom{m}{\ell}a^{m{-}\ell}(\text{ad}\,a)^\ell(c),$$
together with the formula
$$
(\text{ad}\,x^{\beta}\partial')^\ell(x^\alpha\partial)
=x^{\alpha{+}\ell\beta}(a_\ell\partial-b_\ell\partial'),\leqno(16)$$
obtained by induction on $\ell$ when taking $a=x^{\beta}\partial',
\,c=x^\alpha\partial$.
\end{proof}
To simplify formulas in  the sequel, we introduce the operator
$d^{(\ell)}(\ell\geq 0)$ on $U(x^{\eta}\mathbf{S})$ defined by
$d^{(\ell)}:=\frac{1}{\ell!}(\text{\rm ad}\,e)^{\ell}$. From (16)
and the derivation property of $d^{(\ell)}$, it is easy to get
\begin{lemm} For  $\mathbb{Z}^n$-homogeneous elements $x^{\alpha}\partial$, $a_i$, the
following equalities hold in $U(x^{\eta}\mathbf{S}):$
\begin{gather*}
d^{(\ell)}(x^\alpha\partial)=x^{\alpha+\ell\gamma}(A_\ell\partial-B_\ell\partial_0),
\tag{17}\\
d^{(\ell)}(a_1\cdots
a_s)=\sum_{\ell_1{+}\cdots{+}\ell_s=\ell}d^{(\ell_1)}(a_1)\cdots
d^{(\ell_s)}(a_s), \tag{18}
\end{gather*}
where $A_\ell=\frac{1}{\ell!}
\prod\limits_{j=0}^{\ell{-}1}\partial_0(\alpha{+}j\gamma)=\frac{1}{\ell!}
\prod\limits_{j=0}^{\ell{-}1}\partial_0(\alpha{+}j\eta),
\,B_\ell=\partial(\gamma)A_{\ell{-}1}$, and set $A_0=1$, $A_{-1}=0$.
\end{lemm}

Denote by $(U(x^{\eta}\mathbf{S}),\, m,\, \iota,\, \Delta_0,\,
S_0,\, \varepsilon_0)$ the standard Hopf algebra structure of the
universal enveloping algebra $U(x^{\eta}\mathbf{S})$ for the Lie
algebra $x^{\eta}\mathbf{S}$.

\subsection{Quantization of $U(x^{\eta}\mathbf{S})$ in char 0}
We can perform the process of twisting the standard Hopf structure
$(U(x^{\eta}\mathbf{S})[[t]],\, m,\, \iota,\, \Delta_0,\, S_0,\,
\varepsilon_0)$ by the Drinfel'd twist $\mathcal{F}$ constructed in
Proposition 1.9.

\smallskip
The following Lemma is very useful to our main result in this
section.
\begin{lemm}
For $a \in \mathbb{F},\,\alpha \in \mathbb{Z}^n$, and $x^{\a}\pa \in
x^{\a}T_{\a-\eta}$, one has
\begin{gather*}
\bigl((x^{\alpha}\partial)^s\otimes 1\bigr)\cdot
F_a=F_{a+s(\a_{k'}-\al_k)}\cdot
\bigl((x^{\alpha}\partial)^s\otimes 1\bigr), \tag{19}\\
(x^{\alpha}\partial)^s\cdot u_a=u_{a+s(\al_k-\a_{k'})}\cdot
\Bigl(\sum\limits_{\ell=0}^{\infty}d^{(\ell)}\bigl((x^\alpha\partial)^s\bigr)\cdot
h_{1-a}^{\lg \ell\rg}t^\ell\Bigr), \tag{20} \\
\bigl(1\otimes (x^{\alpha}\partial)^s\bigr)\cdot
F_a=\sum\limits_{\ell=0}^{\infty}(-1)^\ell F_{a+\ell}\cdot
\Bigl(h_a^{\lg \ell\rg}\otimes
d^{(\ell)}((x^{\alpha}\partial)^s)t^\ell\Bigr). \tag{21}
\end{gather*}
\end{lemm}
\begin{proof}
For (19): By (14), one has
\begin{equation*}
\begin{split}
(x^\alpha\partial\otimes 1)\cdot F_a
&=\sum\limits_{m=0}^{\infty}\frac{1}{m!}x^\alpha\partial \cdot
h_a^{\lg m\rg}\otimes e^mt^m
\\
&=\sum\limits_{m=0}^{\infty}\frac{1}{m!}h_{a+(\a_{k'}-\al_k)}^{\lg
m\rg}\cdot x^\alpha\partial \otimes e^mt^m
\\
&=F_{a+(\a_{k'}-\al_k)}\cdot(x^\alpha\partial \otimes 1).
\end{split}
\end{equation*}
By induction on $s$, we obtain the result.

For (20):  Let
$a_\ell=\prod\limits_{j=0}^{\ell-1}\partial_0(\alpha{+}j\gamma),\,
b_\ell=\ell\,\partial(\gamma)a_{\ell{-}1}$, using  induction on $s$.
For $s=1$, using (7), (11), (13) and (15), we get
\begin{equation*}
\begin{split}
x^{\alpha}\partial\cdot u_a &=x^{\alpha}\partial\cdot
\left(\sum\limits_{r=0}^{\infty}\frac{(-1)^r}{r!} h_{-a}^{[r]}\cdot
e^rt^r \right) \\
 &=
\sum\limits_{r=0}^{\infty}\frac{(-1)^r}{r!}x^{\alpha}\partial\cdot
h_{-a}^{[r]}\cdot e^rt^r  \\
&=\sum\limits_{r=0}^{\infty}\frac{(-1)^r}{r!}h_{-a-(\al_k-\a_{k'})}^{[r]}\cdot
x^{\alpha}\partial\cdot
 e^rt^r \\ &=\sum\limits_{r=0}^{\infty}\frac{(-1)^r}{r!}h_{-a-(\al_k-\a_{k'})}^{[r]}
 \left(\sum\limits_{\ell=0}^{r}(-1)^\ell\dbinom{r}{\ell} e^{r-\ell}\cdot
 x^{\alpha+\ell\gamma}(a_\ell\partial-b_\ell\partial_0)t^r
 \right)\\
 &=\sum\limits_{r,\ell=0}^{\infty}\frac{({-}1)^{r{+}\ell}}{(r{+}\ell)!}
 h_{{-}a{-}(\al_k-\a_{k'})}^{[r{+}\ell]}
 \left(({-}1)^\ell\dbinom{r{+}\ell}{\ell} e^{r}\cdot
 x^{\alpha{+}\ell\gamma}(a_\ell\partial{-}b_\ell\partial_0)t^{r{+}\ell}
  \right)\\
&=\sum\limits_{r,\ell=0}^{\infty}\frac{(-1)^{r}}{r!\ell!}
  h_{-a-(\al_k-\a_{k'})}^{[r]}\cdot h_{-a-(\al_k-\a_{k'})-r}^{[\ell]}\cdot
 e^{r}\cdot x^{\alpha+\ell\gamma}(a_\ell\partial-b_\ell\partial_0)t^{r+\ell}
  \\
&=\sum\limits_{\ell=0}^{\infty}\left(\sum\limits_{r=0}^{\infty}\frac{(-1)^{r}}{r!}
h_{-a-(\al_k-\a_{k'})}^{[r]} e^{r}t^r\right)
h_{-a-(\al_k-\a_{k'})}^{[\ell]}\cdot
x^{\alpha+\ell\gamma}(A_\ell\partial{-}B_\ell\partial_0)t^{\ell}
  \\
&=u_{a+(\al_k-\a_{k'})}\cdot\sum\limits_{\ell=0}^{\infty}h_{-a-(\al_k-\a_{k'})}^{[\ell]}\cdot
x^{\alpha+\ell\gamma}(A_\ell\partial-B_\ell\partial_0)t^{\ell}
  \\
\end{split}
\end{equation*}
\begin{equation*}
\begin{split}
&=u_{a+(\al_k-\a_{k'})}\cdot\sum\limits_{\ell=0}^{\infty}x^{\alpha+\ell\gamma}
  (A_\ell\partial-B_\ell\partial_0)\cdot h_{-a+\ell}^{[\ell]}t^{\ell}
\\
&=u_{a+(\al_k-\a_{k'})}\cdot\sum\limits_{\ell=0}^{\infty}d^{(\ell)}(x^{\al}\pa)\cdot
h_{-a+1}^{\lg
 \ell\rg}t^{\ell},
\end{split}
\end{equation*}
where $A_\ell=\frac{1}{\ell!}
\prod\limits_{j=0}^{\ell{-}1}\partial_0(\alpha{+}j\gamma)=\frac{1}{\ell!}
\prod\limits_{j=0}^{\ell{-}1}\partial_0(\alpha{+}j\eta),
\,B_\ell=\partial(\gamma)A_{\ell{-}1}$, and set $A_0=1$, $A_{-1}=0$.

Suppose $s\ge 1$. Using Lemma 2.2 and the induction hypothesis on
$s$, we have
\begin{equation*}
\begin{split}
(x^{\alpha}\partial)^{s{+}1}\cdot u_a &= x^{\alpha}\partial\cdot
u_{a+ s(\al_k-\a_{k'})} \cdot
\sum\limits_{n=0}^{\infty}d^{(n)}((x^{\alpha}\partial)^s)\cdot
h_{1-a}^{\langle
n \rangle} t^{n}\\
&=
u_{a{+}(s{+}1)(\al_k-\a_{k'})}{\cdot}\Bigl(\sum\limits_{m=0}^{\infty}d^{(m)}(x^{\alpha}\partial)\cdot
h_{1{-}a{-}s(\al_k-\a_{k'})}^{\langle m \rangle}
t^{m}\Bigr)\\
&\quad\qquad\qquad\qquad\cdot\Bigl(\sum\limits_{n=0}^{\infty}d^{(n)}((x^{\alpha}\partial)^s)\cdot
h_{1{-}a}^{\langle
n \rangle} t^{n}\Bigr)\\
&=u_{a{+}(s{+}1)(\al_k-\a_{k'})}\cdot\Bigl(\sum\limits_{m,n=0}^{\infty}d^{(m)}(x^{\alpha}\partial)
d^{(n)}((x^{\alpha}\partial)^s)h_{1-a+n}^{\langle m
\rangle} h_{1-a}^{\langle n \rangle} t^{n+m}\Bigr)\\
&=u_{a{+}(s{+}1)(\al_k-\a_{k'})}\cdot\Bigl(\sum\limits_{\ell=0}^{\infty}\sum\limits_{m+n=\ell}d^{(m)}(x^{\alpha}\partial)
d^{(n)}((x^{\alpha}\partial)^s) h_{1-a}^{\langle \ell \rangle} t^{\ell}\Bigr)\\
&=u_{a{+}(s{+}1)(\al_k-\a_{k'})}\cdot\Bigl(\sum\limits_{\ell=0}^{\infty}d^{(\ell)}((x^{\alpha}\partial)^{s{+}1})h_{1-a}^{\langle
\ell \rangle} t^{\ell}\Bigr),
\end{split}
\end{equation*}
where we get the first and second  ``=" by using the inductive
hypothesis, the third by using (14) \& (18) and the fourth by using
(6) \& (18).

For (21): For $s$=1, using (15) we get
\begin{equation*}
\begin{split}
(1 \otimes x^\alpha\partial )\cdot F_a
&=\sum\limits_{m=0}^{\infty}\frac{1}{m!} h_a^{\lg m\rg}\otimes
x^\alpha\partial \cdot e^mt^m
\\
&=\sum\limits_{m=0}^{\infty}\frac{1}{m!}h_{a}^{\lg m\rg}\otimes
\left( \sum\limits_{\ell=0}^{m}(-1)^\ell\dbinom{m}{\ell}
e^{m{-}\ell}\cdot
x^{\alpha{+}\ell\gamma}(a_\ell\partial{-}b_\ell\partial_0)t^m\right)\\
&=\sum\limits_{m=0}^{\infty}\sum\limits_{\ell=0}^{\infty}(-1)^\ell\frac{1}{m!\ell!}h_{a}^{\lg
m+\ell\rg}\otimes e^{m}\cdot
x^{\alpha+\ell\gamma}(a_\ell\partial-b_\ell\partial_0)t^{m+\ell}
\\
&=\sum\limits_{\ell=0}^{\infty}(-1)^\ell\left(\sum\limits_{m=0}^{\infty}\frac{1}{m!}h_{a+\ell}^{\lg
m\rg}\otimes e^mt^m\right)\Bigl(h_a^{\lg \ell\rg}\otimes
d^{(\ell)}(x^{\al}\pa) t^{\ell}\Bigr)
\\
&=\sum\limits_{\ell=0}^{\infty}(-1)^\ell F_{a+\ell}\cdot
\Bigl(h_a^{\lg \ell\rg}\otimes d^{(\ell)}(x^{\al}\pa)t^{\ell}\Bigr).
\end{split}
\end{equation*}
For $s>1$, it follows from the induction hypothesis  \& (18).
\end{proof}

The following theorem gives the quantization of
$U(x^{\eta}\mathbf{S})$ by Drinfel'd twist $\mathcal{F}$, which is
essentially determined by the Lie bialgebra triangular structure on
$x^{\eta}\mathbf{S}$.
\begin{theorem}
Fix two distinguished elements $h=\partial_k{-}\pa_{k'}$,
$e=x^{\gamma}\partial_0$, where
 $\gamma$ satisfies $\gamma_k{-}\gamma_{k'}=1$ such that
$[h,\,e]=e$
 in the generalized Cartan type $\mathbf{S}$ Lie
algebra $x^{\eta}\mathbf{S}$ over $\mathbb{F}$, there exists a
structure of noncommutative and noncocommutative Hopf algebra
$(U(x^{\eta}\mathbf{S})[[t]],m,\iota,\Delta,S,\varepsilon)$ on
$U(x^{\eta}\mathbf{S})[[t]]$ over $\mathbb{F}[[t]]$
 with
$U(x^{\eta}\mathbf{S})[[t]]/tU(x^{\eta}\mathbf{S})[[t]]$ $\cong
U(x^{\eta}\mathbf{S})$, which leaves the product of
$U(x^{\eta}\mathbf{S})[[t]]$ undeformed but with the deformed
coproduct, antipode and counit defined by
\begin{gather*}
\Delta(x^{\alpha}\partial)=x^{\alpha}\partial\otimes
(1{-}et)^{\al_k-\a_{k'}}+\sum\limits_{\ell=0}^{\infty}(-1)^\ell
h^{\lg \ell\rg}\otimes (1{-}et)^{-\ell}\cdot
d^{(\ell)}(x^{\al}\pa)t^\ell,
\tag{22}\\
S(x^{\alpha}\partial)=-(1{-}et)^{-(\al_k-\a_{k'})}\cdot\Bigl(\sum\limits_{\ell=0}^{\infty}
d^{(\ell)}(x^{\al}\pa)\cdot h_1^{\lg \ell\rg}t^\ell\Bigr),
\tag{23}\\
\varepsilon(x^{\alpha}\partial)=0,  \tag{24}
\end{gather*}
where $x^{\alpha}\partial \in x^{\al}T_{\al-\eta}$.
\end{theorem}
\begin{proof}
By Lemmas 1.4 and 1.6, it follows from (19) and (21) that
\begin{equation*}
\begin{split}
\Delta(x^{\alpha}\partial)
&=\mathcal{F}\cdot\Delta_0(x^{\alpha}\partial)\cdot\mathcal{F}^{-1}\\
&=\mathcal{F}\cdot(x^{\alpha}\partial\otimes 1)\cdot F+
\mathcal{F}\cdot(1\otimes x^{\alpha}\partial)\cdot F
\\
&=\Bigl(\mathcal{F}
F_{\al_{k'}-\a_k}\Bigr)\cdot(x^{\alpha}\partial{\otimes} 1)+
\sum\limits_{\ell=0}^{\infty}(-1)^\ell
\Bigl(\mathcal{F}F_{\ell}\Bigr)\cdot\Bigl(h^{\lg \ell\rg}{\otimes}
d^{(\ell)}(x^{\al}\pa)t^\ell\Bigr)
\\ &=\Bigl(1\otimes
(1{-}et)^{\al_k-\a_{k'}}\Bigr)\cdot (x^{\alpha}\partial\otimes
1)\\
&\quad+\sum\limits_{\ell=0}^{\infty}(-1)^\ell \Bigl(1\otimes
(1{-}et)^{-\ell}\Bigr)\cdot\Bigl( h^{\lg \ell\rg}\otimes
d^{(\ell)}(x^{\al}\pa)t^\ell\Bigr)
\\ &=x^{\alpha}\partial\otimes
(1{-}et)^{\al_k-\a_{k'}}+\sum\limits_{\ell=0}^{\infty}(-1)^\ell
h^{\lg \ell\rg}\otimes (1{-}et)^{-\ell}\cdot
d^{(\ell)}(x^{\al}\pa)t^\ell.
\\
\end{split}
\end{equation*}
By (20) and Lemma 1.6, we obtain
\begin{equation*}
\begin{split}
S(x^{\alpha}\partial)&=u^{-1}S_0(x^{\alpha}\partial)\,u=-v\cdot
x^{\alpha}\partial\cdot u\\ &=-v \cdot u_{\al_k-\a_{k'}}\cdot
\Bigl(\sum\limits_{\ell=0}^{\infty}d^{(\ell)}(x^{\al}\pa)\cdot
h_{1}^{\lg \ell\rg}t^\ell\Bigr)
\\ &=-(1{-}et)^{-(\al_k-\a_{k'})}\cdot\Bigl(\sum\limits_{\ell=0}^{\infty}
d^{(\ell)}(x^{\al}\pa)\cdot h_1^{\lg \ell\rg}t^\ell\Bigr).
\end{split}
\end{equation*}

Hence, we get the result.
\end{proof}

For later use, we need to make the following
\begin{lemm} For $s\ge 1$, one has
\begin{gather*}
\Delta((x^\alpha\partial)^s)=\sum_{0\le j\le s\atop
\ell\ge0}\dbinom{s}{j}({-}1)^\ell(x^{\alpha}\partial)^jh^{\langle
\ell\rangle}\otimes(1{-}et)^{j(\al_k-\a_{k'}){-}\ell}
d^{(\ell)}((x^\al\partial)^{s{-}j})t^\ell.\tag{\text{\rm i}}\\
S((x^{\alpha}\partial)^s)=
(-1)^s(1{-}et)^{-s(\al_k-\a_{k'})}\cdot\Bigl(\sum\limits_{\ell=0}^{\infty}
d^{(\ell)}((x^{\alpha}\partial)^s)\cdot h_1^{\lg
\ell\rg}t^\ell\Bigr).\tag{\text{\rm ii}}
\end{gather*}
\end{lemm}
\begin{proof} By (19), (21) and Lemma 1.6, we obtain
\begin{equation*}
\begin{split}
\Delta((x^{\alpha}\partial)^s)&=\mathcal{F}\Bigl(x^{\al}\partial\otimes
1+1\otimes x^{\alpha}\partial\Bigr)^s\mathcal{F}^{-1}\\
&=\sum_{j=0}^s\binom{s}{j}\mathcal{F}F_{j(\al_{k'}{-}\a_k)}
(x^\alpha\partial{\otimes}
1)^j\Bigl(\sum_{\ell\ge0}({-}1)^\ell\mathcal{F}F_{\ell}\bigl(h^{\langle\ell\rangle}{\otimes}
d^{(\ell)}((x^\alpha\partial)^{s{-}j})t^\ell\bigr)\Bigr)\\
&=\sum_{j=0}^s\sum_{\ell\ge0}\binom{s}{j}({-}1)^\ell\bigl((x^\alpha\partial)^j{\otimes}
(1{-}et)^{j(\al_k{-}\a_{k'}){-}\ell}\bigr)\bigl(h^{\langle\ell\rangle}{\otimes}
d^{(\ell)}((x^\alpha\partial)^{s{-}j})t^\ell\bigr)\\
&=\sum_{0\le j\le s\atop
\ell\ge0}\dbinom{s}{j}({-}1)^\ell(x^{\alpha}\partial)^jh^{\langle
\ell\rangle}\otimes(1{-}et)^{j(\al_k-\a_{k'}){-}\ell}
d^{(\ell)}((x^\al\partial)^{s{-}j})t^\ell.
\end{split}
\end{equation*}
Again by (20) and Lemma 1.6, we get
\begin{equation*}
\begin{split}
S((x^{\alpha}\partial)^s)&=u^{-1}S_0((x^{\alpha}\partial)^s)\,u=(-1)^s
v\cdot (x^{\alpha}\partial)^s\cdot u\\ &=(-1)^s v \cdot
u_{s(\al_k-\a_{k'})}\cdot
\Bigl(\sum\limits_{\ell=0}^{\infty}d^{(\ell)}((x^{\alpha}\partial)^s)\cdot
h_{1}^{\lg \ell\rg}t^\ell\Bigr)
\\ &=(-1)^s(1{-}et)^{-s(\al_k-\a_{k'})}\cdot\Bigl(\sum\limits_{\ell=0}^{\infty}
d^{(\ell)}((x^{\alpha}\partial)^s)\cdot h_1^{\lg
\ell\rg}t^\ell\Bigr).
\end{split}
\end{equation*}

So this completes the proof.
\end{proof}

\subsection{Quantization integral forms of $\mathbb{Z}$-form $\mathbf{S}_{\mathbb{Z}}^+$ in
char $0$} As we known,
$\{\al_nx^{\al-\ep_n}D_i-\al_ix^{\al-\ep_i}D_n=x^{\al-\ep_i-\ep_n}(\al_n\pa_i-\al_i\pa_n)
\mid \alpha \in\mathbb{Z}_+^n, 1\le i< n\,\}$ is a
$\mathbb{Z}$-basis of $\mathbf{S}_{\mathbb{Z}}^+$, as a subalgebra
of both the simple Lie $\mathbb{Z}$-algebras
$x^{-\underline{1}}\mathbf{S}_{\mathbb{Z}}$ and
$\mathbf{W}^+_{\mathbb{Z}}$. In order to get the quantization
integral forms of $\mathbb{Z}$-form $\mathbf{S}_{\mathbb{Z}}^+$, it
suffices to consider what conditions are needed for those
coefficients occurred in the formulae (22) \& (23) to be integral
for the indicated basis elements.

\begin{lemm} $($\cite{CG}$)$ \
For any $a,\,k,\,\ell\in\mathbb{Z}$,
$a^\ell\prod\limits_{j=0}^{\ell-1}(k{+}ja)/\ell!$ is an
integer.\hfill\qed
\end{lemm}
From this Lemma (due to Grunspan), we see that if we take
$\partial_0(\ga)=\pm1$, then $A_\ell$ and $B_\ell$ are integers in
Theorem 2.4. In this paper, the cases we are interested in are:
$\mathrm{(i)}\ h=\partial_k{-}\pa_{k'}$,
$e=x^{\epsilon_k}(\partial_k{-}2\pa_{k'})$ $(1\leq k \neq k'\leq
n)$; $\mathrm{(ii)}\ h=\partial_k{-}\pa_{k'}$,
$e=x^{\epsilon_k-\epsilon_m}\partial_m$ $(1\leq k \neq k'\neq m\leq
n)$. The latter will be discussed in Section 5. Denote by
$\mathcal{F}(k,k')$ the corresponding Drinfel'd twist in the case
$\mathrm{(i)}$. As a result of Theorem 2.4, we have
\begin{coro}
Fix distinguished elements $h:=\pa_k{-}\pa_{k'}$,
$e:=x^{\epsilon_k}(\pa_k{-}2\pa_{k'})$ $(1\leq k \neq k'\leq n)$,
the corresponding quantization of $U(\mathbf{S}^+_{\mathbb{Z}})$
over $U(\mathbf{S}^+_{\mathbb{Z}})[[t]]$ by Drinfel'd twist
$\mathcal{F}(k,k')$ with the product undeformed is given by
\begin{gather*}
\Delta(x^{\alpha}\pa)=x^{\alpha}\pa\otimes
(1{-}et)^{\alpha_k{-}\a_{k'}}+\sum\limits_{\ell=0}^{\infty}{({-}1)}^\ell
\,h^{\lg \ell\rg}\otimes (1{-}et)^{{-}\ell}\cdot\tag{25}  \\
\qquad\qquad\qquad\qquad\qquad\qquad\qquad\qquad\qquad\, \cdot\,
x^{\alpha{+}\ell\epsilon_k}\bigl(A_\ell\pa-B_\ell(\pa_k{-}2\pa_{k'})\bigr)t^\ell,\\
S(x^{\alpha}\pa)={-}(1{-}et)^{-(\alpha_k{-}\a_{k'})}\cdot\Bigl(\sum\limits_{\ell=0}^{\infty}
\,x^{\alpha{+}\ell\epsilon_k}\bigl(A_\ell\pa-B_\ell(\pa_k{-}2\pa_{k'})\bigr)\cdot
h_1^{\lg \ell\rg}t^\ell\Bigr), \tag{26}
\end{gather*}
\begin{gather*}
\varepsilon(x^{\alpha}\pa)=0,  \tag{27}
\end{gather*}
where $
A_\ell=\frac{1}{\ell!}\prod\limits_{j=0}^{\ell-1}(\al_k{-}2\a_{k'}{+}j),\,
B_\ell=\pa(\epsilon_k) A_{\ell{-}1}$ with $A_0=1, A_{-1}=0$.
\end{coro}
\begin{remark} We get $n(n-1)$ {\it basic
Drinfel'd twists} $\mathcal{F}(1,2),{\cdots},\mathcal{F}(1,n)$,
 $\mathcal{F}(2,1),$  $\cdots,\mathcal{F}(n,n-1)$ over $U(\mathbf{S}^+_{\mathbb{Z}})$. It is interesting to
consider the products of some {\it basic Drinfel'd twists}, using
the same argument as the proof of Theorem 2.4, one can get many more
new Drinfel'd twists (which depends on a bit more calculations to be
done), which will lead to many more new complicated quantizations
not only over the $U(\mathbf{S}_{\mathbb{Z}}^+)[[t]]$, but the
possible quantizations over the
$\mathbf{u}_{t,q}(\mathbf{S}(n;\underline{1}))$ as well, via our
modulo reduction approach developed in the next section.
\end{remark}

\section{Quantizations of the special algebra $\mathbf{S}(n;\underline{1})$ }

In this section, firstly, we make {\it modulo $p$ reduction and base
change with the $\mathcal K[[t]]$ replaced by $\mathcal K[t]$}, for
the quantization of $U(\mathbf{S}^+_{\mathbb{Z}})$ in char $0$
(Corollary 2.7) to yield the quantization of
$U(\mathbf{S}(n;\underline{1}))$, for the restricted simple modular
Lie algebra $\mathbf{S}(n;\underline{1})$ in char $p$. Secondly, we
shall further make {\it ``$p$-restrictedness" reduction as well as
base change with the $\mathcal K[t]$ replaced by $\mathcal
K[t]_p^{(q)}$}, for the quantization of
$U(\mathbf{S}(n;\underline{1}))$, which will lead to the required
quantization of $\mathbf{u}(\mathbf{S}(n;\underline{1}))$, the
restricted universal enveloping algebra of
$\mathbf{S}(n;\underline{1})$.

\subsection{Modulo $p\,$ reduction and base change}
Let $\mathbb{Z}_p$ be the prime subfield of $\mathcal{K}$ with
$\text{char}(\mathcal{K})=p$. When considering
$\mathbf{W}_{\mathbb{Z}}^+$ as a $\mathbb{Z}_p$-Lie algebra, namely,
making modulo $p$ reduction for the defining relations of
$\mathbf{W}_{\mathbb{Z}}^+$, denoted by
$\mathbf{W}_{\mathbb{Z}_p}^+$, we see that
$(J_{\underline{1}})_{\mathbb{Z}_p}=\text{Span}_{\mathbb{Z}_p}\{x^\alpha
D_i \mid \exists\, j: \alpha_j\ge p\,\}$ is a maximal ideal of
$\mathbf{W}^+_{\mathbb{Z}_p}$, and
$\mathbf{W}^+_{\mathbb{Z}_p}/(J_{\underline{1}})_{\mathbb{Z}_p}
\cong \mathbf{W}(n;\underline{1})_{\mathbb{Z}_p}
=\text{Span}_{\mathbb{Z}_p}\{x^{(\alpha)}D_{i}\mid 0\le \alpha\le
\tau, 1\le i\le n\}$. For the subalgebra
$\mathbf{S}_{\mathbb{Z}}^+$, we have
$\mathbf{S}^+_{\mathbb{Z}_p}/(\mathbf{S}^+_{\mathbb{Z}_p}\cap(J_{\underline{1}})_{\mathbb{Z}_p})
\cong \mathbf{S}'(n;\underline{1})_{\mathbb{Z}_p}$. We denote simply
$\mathbf{S}^+_{\mathbb{Z}_p}\cap(J_{\underline{1}})_{\mathbb{Z}_p}$
as  $(J^+_{\underline{1}})_{\mathbb{Z}_p}$.

Moreover, we have $\mathbf{S}'(n;\underline{1})
=\mathcal{K}\otimes_{\mathbb{Z}_p}\mathbf{S}'(n;\underline{1})_{\mathbb{Z}_p}
=\mathcal{K}\mathbf{S}'(n;\underline{1})_{\mathbb{Z}_p}$, and
$\mathbf{S}^+_{\mathcal{K}}=\mathcal{K}\mathbf{S}^+_{\mathbb{Z}_p}$.

Observe that the ideal
$J^+_{\underline{1}}:=\mathcal{K}(J^+_{\underline{1}})_{\mathbb{Z}_p}$
generates an ideal of $U(\mathbf{S}^+_{\mathcal{K}})$ over
$\mathcal{K}$, denoted by
$J:=J^+_{\underline{1}}U(\mathbf{S}^+_{\mathcal{K}})$, where
$\mathbf{S}^+_{\mathcal{K}}/J^+_{\underline{1}}\cong
\mathbf{S}'(n;\underline{1})$. Based on the formulae (25) \& (26),
$J$ is a Hopf ideal of $U(\mathbf{S}^+_{\mathcal{K}})$ satisfying
$U(\mathbf{S}^+_{\mathcal{K}})/J\cong
U(\mathbf{S}'(n;\underline{1}))$. Note that elements $\sum a_{i,
\alpha}\frac{1}{\alpha!}x^{\alpha}D_i$ in
$\mathbf{S}^+_{\mathcal{K}}$ for $0\le\alpha\le\tau$ will be
identified with $\sum a_{i, \alpha}x^{(\alpha)}D_i$ in
$\mathbf{S}'(n;\underline{1})$ and those in $J_{\underline{1}}$ with
$0$. Hence, by Lemma 1.2 and Corollary 2.7, we get the quantization
of $U(\mathbf{S}'(n;\underline{1}))$ over
$U_t(\mathbf{S}'(n;\underline{1})):=U(\mathbf{S}'(n;\underline{1}))\otimes_{\mathcal
K}\mathcal K[t]$ (not necessarily in
$U(\mathbf{S}'(n;\underline{1}))[[t]]$, as seen in formulae (28) \&
(29)) as follows.
\begin{theorem}
Fix two distinguished elements
$h:=D_{kk'}(x^{(\epsilon_k+\epsilon_{k'})})$,
$e:=2D_{kk'}(x^{(2\epsilon_k+\epsilon_{k'})})$ $(1\leq k\neq k'\leq
n)$, the corresponding quantization of
$U(\mathbf{S}'(n;\underline{1}))$ over
$U_t(\mathbf{S}'(n;\underline{1}))$ with the product undeformed is
given by
\begin{gather*}
\Delta(D_{ij}(x^{(\alpha)}))=D_{ij}(x^{(\alpha)})\otimes
(1{-}et)^{\alpha_k{-}\delta_{ik}{-}\delta_{jk}-\a_{k'}{+}\delta_{ik'}{+}\delta_{jk'}}\tag{28}\\
\qquad\qquad\qquad\quad\,+\sum\limits_{\ell=0}^{p{-}1}{({-}1)}^\ell
h^{\lg \ell\rg}\otimes(1{-}et)^{{-}\ell}\Bigl(\bar{A}_\ell
D_{ij}(x^{(\alpha{+}\ell\epsilon_k)})\\
\qquad\qquad\qquad\qquad\qquad\qquad\qquad\, +\,\bar{B}_\ell
(\delta_{ik}D_{k'j}{+}\delta_{jk}D_{ik'})
(x^{(\alpha{+}(\ell{-}1)\epsilon_k+\epsilon_{k'})})\Bigr)t^\ell,
\end{gather*}
\begin{gather*}
S(D_{ij}(x^{(\alpha)}))={-}(1{-}et)^{-\alpha_k{+}\delta_{ik}{+}\delta_{jk}+\a_{k'}{-}\delta_{ik'}{-}\delta_{jk'}}
\cdot\Bigl(\sum\limits_{\ell=0}^{p{-}1}\bigl(\bar{A}_\ell
D_{ij}(x^{(\alpha{+}\ell\epsilon_k)})   \tag{29}\\
\qquad\qquad +\bar{B}_\ell
(\delta_{ik}D_{k'j}{+}\delta_{jk}D_{ik'})(x^{(\alpha{+}(\ell-1)\epsilon_k+\epsilon_{k'})})\bigr)
\cdot h_1^{\lg \ell\rg}t^\ell\Bigr),
\\
\varepsilon(D_{ij}(x^{(\alpha)}))=0, \tag{30}\\
\Delta(x^{(\tau-(p-1)\epsilon_j)}D_j)=x^{(\tau-(p-1)\epsilon_j)}D_j{\otimes}
(1{-}et)^{p\,(\delta_{jk'}{-}\delta_{jk})}+1{\otimes}
x^{(\tau-(p-1)\epsilon_j)}D_j, \tag {31}\\
S(x^{(\tau-(p-1)\epsilon_j)}D_j)=-(1{-}et)^{p\,(\delta_{jk}{-}\delta_{jk'})}x^{(\tau-(p-1)\epsilon_j)}D_j,
\tag {32}\\
\varepsilon(x^{(\tau-(p-1)\epsilon_j)}D_j)=0, \tag {33}
\end{gather*}
where $0\le \alpha \le \tau$, $1\le j<i\le n$,  $\bar
A_\ell=\ell!\binom{\alpha_k{+}\ell}{\ell}(A_\ell-\delta_{jk}A_{\ell-1}-\delta_{ik}A_{\ell-1})\,(\text{\rm
mod} \,p)$, $\bar B_\ell=2\ell!
\binom{\alpha_k{+}\ell-1}{\ell-1}(\a_{k'}+1)A_{\ell-1}\,(\text{\rm
mod}\,p)$,
$A_\ell=\frac{1}{\ell!}\prod\limits_{m=0}^{\ell-1}(\al_k-\delta_{jk}-\delta_{ik}
-2\a_{k'}+2\delta_{jk'}+2\delta_{ik'}+m)$ and $A_0=1, A_{-1}=0$.
\end{theorem}

Note that (28), (29) \& (30) give the corresponding quantization of
$U(\mathbf{S}(n;\underline{1}))$ over
$U_t(\mathbf{S}(n;\underline{1})):=U(\mathbf{S}(n;\underline{1}))\otimes_{\mathcal
K}\mathcal K[t]$ (also over $U(\mathbf{S}(n;\underline{1}))[[t]]$).
It should be noticed that in this step --- inducing from the
quantization integral form of $U(\mathbf S_{\mathbb Z}^+)$ and
making the modulo $p$ reduction, we used the first base change with
$\mathcal K[[t]]$ replaced by $\mathcal K[t]$, and the objects from
$U(\mathbf S(n;\underline1))[[t]]$ turning to $U_t(\mathbf
S(n;\underline1))$.

\subsection{Modulo ``$p$-restrictedness" reduction and base change}
Let $I$ be the ideal of $U(\mathbf{S}(n;\underline{1}))$ over
$\mathcal{K}$ generated by
$(D_{ij}(x^{(\epsilon_i+\epsilon_j)}))^p-D_{ij}(x^{(\epsilon_i+\epsilon_j)})$
and $(D_{ij}(x^{(\al)}))^p$ with $\alpha\ne \epsilon_i+\epsilon_j$
for $0\le\alpha\le \tau$ and $1\le j<i\le n$.
$\mathbf{u}(\mathbf{S}(n;\underline{1}))=U(\mathbf{S}(n;\underline{1}))/I$
is of dimension $p^{(n-1)(p^n-1)}$. In order to get a reasonable
quantization of finite dimension for
$\mathbf{u}(\mathbf{S}(n;\underline{1}))$ in char $p$, at first, it
is necessary to clarify in concept the underlying vector space in
which the required $t$-deformed object exists. According to our
modular reduction approach, it should start to be induced from the
$\mathcal{K}[t]$-algebra $U_t(\mathbf{S}(n;\underline{1}))$ in
Theorem 3.1.

Firstly, we observe the following fact
\begin{lemm} $(\text{\rm i})$ \ $(1-et)^p\equiv 1 \quad (\text{\rm mod}\,p, I)$.

\smallskip $(\text{\rm ii})$ \ $(1-et)^{-1}\equiv
1+et+\cdots+e^{p-1}t^{p-1} \quad (\text{\rm mod}\,
 p, I)$.

\smallskip $(\text{\rm iii})$ \ $h_a^{\lg \ell\rg} \equiv 0 \quad
(\text{\rm mod} \, p, I)$ for $\ell \geq p$, and $a\in\mathbb{Z}_p$.
\end{lemm}
\begin{proof} (i), (ii) follow from $e^p=0$ in $\mathbf{u}(\mathbf{S}(n;\underline{1}))$.

(iii)  For $\ell\in\mathbb{Z}_+$, there is a unique decomposition
$\ell=\ell_0+\ell_1p$ with $0\le \ell_0<p$ and $\ell_1\ge 0$. Using
the formulae (4) \& (6), we have
\begin{equation*}
\begin{split}
h_a^{\langle \ell\rangle}& =h_a^{\langle \ell_0\rangle}\cdot
h_{a+\ell_0}^{\langle \ell_1p\rangle}\equiv h_a^{\langle
\ell_0\rangle}\cdot (h_{a+\ell_0}^{\langle
p\rangle})^{\ell_1}\,\qquad
(\text{mod } p)\\
&\equiv h_a^{\langle \ell_0\rangle}\cdot (h^p-h)^{\ell_1}\quad
(\text{mod } p),
\end{split}
\end{equation*}
where we used the facts that $(x+1)(x+2)\cdots(x+p-1)\equiv
x^{p-1}-1\; (\text{mod } p)$, and $(x+a+\ell_0)^p\equiv
x^p+a+\ell_0\ (\text{mod } p)$. Hence, $h_a^{\lg \ell\rg} \equiv
0$ (mod $p, \,I$) for $\ell \geq p$.
\end{proof}

The above Lemma, together with Theorem 3.1, indicates that the
required $t$-deformation of
$\mathbf{u}(\mathbf{S}(n;\underline{1}))$ (if it exists) in fact
only happens in a $p$-truncated polynomial ring (with degrees of $t$
less than $p$) with coefficients in
$\mathbf{u}(\mathbf{S}(n;\underline{1}))$, i.e.,
$\mathbf{u}_{t,q}(\mathbf{S}(n;\underline{1})):=
\mathbf{u}(\mathbf{S}(n;\underline{1}))\otimes_{\mathcal K}
\mathcal{K}[t]_p^{(q)}$ (rather than in
$\mathbf{u}_t(\mathbf{S}(n;\underline{1})):=\mathbf{u}(\mathbf{S}(n;\underline{1}))
\otimes_{\mathcal K}\mathcal K[t]$), where $\mathcal{K}[t]_p^{(q)}$
is taken to be a $p$-truncated polynomial ring which is a quotient
of $\mathcal K[t]$ defined as
$$
\mathcal{K}[t]_p^{(q)}= \mathcal{K}[t]/(t^p-qt), \qquad\text{\it for
}\ q\in\mathcal{K}.\leqno(34)
$$
Thereby, we obtain the underlying ring for our required
$t$-deformation of $\mathbf{u}(\mathbf{S}(n;\underline{1}))$ over
$\mathcal{K}[t]_p^{(q)}$, and
$\dim_{\mathcal{K}}\mathbf{u}_{t,q}(\mathbf{S}(n;\underline{1}))
=p\cdot\dim_{\mathcal{K}}\mathbf{u}(\mathbf{S}(n;\underline{1}))
=p^{1+(n-1)(p^n-1)}$. Via modulo ``restrictedness" reduction, it is
necessary for us to work over the objects from $U_t(\mathbf
S(n;\underline1))$ passage to $U_{t,q}(\mathbf S(n;\underline1))$
first, and then to $\mathbf u_{t,q}(\mathbf S(n;\underline1))$ (see
the proof of Theorem 3.5 below), here we used the second base change
with $\mathcal K[t]_p^{(q)}$ instead of $\mathcal K[t]$.

\smallskip
We are now in a position to describe the following

\begin{defi} With notations as above. A Hopf algebra
$(\mathbf{u}_{t,q}(\mathbf{S}(n;\underline{1}))$, $m,
\iota,\Delta,S,\varepsilon)$ over a ring $\mathcal K[t]_p^{(q)}$ of
characteristic $p$ is said to be a finite-dimensional quantization
of $\mathbf{u}(\mathbf{S}(n;\underline{1}))$ if its Hopf algebra
structure, via modular reduction and base changes, inherits from a
twisting of the standard Hopf algebra $U(\mathbf {S}^+_\mathbb
Z)[[t]]$ by a Drinfeld twist such that
$\mathbf{u}_{t,q}(\mathbf{S}(n;\underline{1}))/t\mathbf{u}_{t,q}(\mathbf{S}(n;\underline{1}))
$ $\cong \mathbf{u}(\mathbf{S}(n;\underline{1}))$.
\end{defi}

To describe $\mathbf{u}_{t,q}(\mathbf{S}(n;\underline{1}))$
explicitly, we still need an auxiliary Lemma.

\begin{lemm} Let $e=2D_{kk'}(x^{(2\epsilon_k+\epsilon_{k'})})$ and $d^{(\ell)}=\frac{1}{\ell!}(\text{\rm
ad}\,e)^\ell$. Then

\smallskip
$(\text{\rm i})$ \ $d^{(\ell)}(D_{ij}(x^{(\alpha)}))=\bar{A}_\ell
D_{ij}(x^{(\alpha{+}\ell\epsilon_k)})+\bar{B}_\ell
(\delta_{ik}D_{k'j}+\delta_{jk}D_{ik'})(x^{(\alpha{+}(\ell-1)\epsilon_k+\epsilon_{k'})})$,

\smallskip\hskip0.5cm
where $\bar A_\ell, \bar B_\ell$ as in Theorem 3.1.

\smallskip $(\text{\rm ii})$ \
$d^{(\ell)}(D_{ij}(x^{(\epsilon_i+\epsilon_j)}))=\delta_{\ell,0}D_{ij}(x^{(\epsilon_i+\epsilon_j)})
-\delta_{1,\ell}(\delta_{ik}-\delta_{jk})e$.

\smallskip $(\text{\rm iii})$ \
$d^{(\ell)}((D_{ij}(x^{(\alpha)}))^p)=\delta_{\ell,0}(D_{ij}(x^{(\alpha)}))^p-\delta_{1,\ell}
(\delta_{ik}-\delta_{jk})\delta_{\alpha,\epsilon_i+\epsilon_j}e$.
\end{lemm}
\begin{proof} (i) Note that
$A_\ell=\frac{1}{\ell!}\prod\limits_{m=0}^{\ell-1}(\al_k{-}\delta_{jk}{-}\delta_{ik}
{-}2\a_{k'}{+}2\delta_{jk'}{+}2\delta_{ik'}{+}m)$, for $0\le
\alpha\le\tau$. By (17) and Theorem 3.1, we get
\begin{equation*}
\begin{split}
d^{(\ell)}(D_{ij}(x^{(\alpha)}))&=\frac{1}{\alpha!}d^{(\ell)}
(x^{\alpha{-}\epsilon_i{-}\epsilon_j}(\al_j\partial_i
{-}\al_i\pa_j))\\
&=\frac{1}{\alpha!}x^{\alpha{-}\epsilon_i{-}\epsilon_j+\ell\epsilon_k}(A_\ell(\al_j\partial_i
{-}\al_i\pa_j)-(\al_j\delta_{ik}{-}\al_i\delta_{jk})A_{\ell{-}1}(\partial_k{-}2\pa_{k'}))\\
&=\bar{A}_\ell D_{ij}(x^{(\alpha{+}\ell\epsilon_k)})+\bar{B}_\ell
(\delta_{ik}D_{k'j}{+}\delta_{jk}D_{ik'})(x^{(\alpha{+}(\ell-1)\epsilon_k+\epsilon_{k'})}).
\end{split}
\end{equation*}

(ii) Note that $A_0=1$ and $A_\ell=0 $ for $\ell \geq 1$,
\begin{gather*}
\bar
A_\ell=\ell!\binom{\alpha_k{+}\ell}{\ell}(A_\ell-\delta_{jk}A_{\ell-1}-\delta_{ik}A_{\ell-1})\quad
(\text{\rm
mod} \,p), \\
\bar B_\ell=2\ell!
\binom{\alpha_k{+}\ell{-}1}{\ell{-}1}(\a_{k'}{+}1)A_{\ell-1}\quad
(\text{\rm mod}\,p).
\end{gather*}
We obtain $\bar A_0=1$ and $\bar B_0=0$. We also obtain $\bar
A_1=-(\delta_{ik}+\delta_{jk})(\al_k+1)$, $\bar B_1=2(\a_{k'}+1)$
and $\bar A_\ell=\bar B_\ell=0$ for $\ell\ge 2$, namely,
$d^{(\ell)}(D_{ij}(x^{(\epsilon_i+\epsilon_j)}))=0$ for $\ell\ge 2$.
So by (i), we have
\begin{equation*}
\begin{split}
d^{(1)}(D_{ij}(x^{(\epsilon_i+\epsilon_j)}))&=-(\delta_{ik}+\delta_{jk})(\al_k+1)
D_{ij}(x^{(\epsilon_i+\epsilon_j+\epsilon_k)})\\
& \quad +2(\a_{k'}+1)(\delta_{ik}D_{k'j}+\delta_{jk}D_{ik'})
(x^{(\epsilon_i+\epsilon_j+\epsilon_{k'})})
\\&=-(\delta_{ik}-\delta_{jk})e.
\end{split}
\end{equation*}
In any case, we arrive at the result as required.

\smallskip
(iii) From (15), we obtain that for $0\le \alpha\le\tau$,
\begin{equation*}
\begin{split}
d^{(1)}\,((D_{ij}(x^{(\alpha)}))^p)&=\frac{1}{(\alpha!)^p}\bigl[\,e,(D_{ij}(x^{\alpha}))^p\,\bigr]
=\frac{1}{(\alpha!)^p}\bigl[\,e,(x^{\alpha-\epsilon_i-\epsilon_j}(\al_j\partial_i-\al_i\pa_j))^p\,\bigr]\\
&=\frac{1}{(\alpha!)^p}\sum\limits_{\ell=1}^p(-1)^\ell\dbinom{p}
{\ell}(x^{\alpha-\epsilon_i-\epsilon_j}(\al_j\partial_i-\al_i\pa_j))^{p-\ell} \\
& \quad\quad \cdot
x^{\epsilon_k{+}\ell(\alpha{-}\epsilon_i-\epsilon_j)} (a_\ell
(\partial_k-2\pa_{k'})-b_\ell(\al_j\partial_i-\al_i\pa_j))\\
&\equiv
-\frac{a_p}{\alpha!}\,x^{2\epsilon_k{+}p(\alpha{-}\epsilon_i-\epsilon_j)}
(\pa_k-2\pa_{k'})\qquad(\text{mod }\,p\,)\\
&\equiv \begin{cases} -{a_p}\,e,\qquad & \text{\it if }\quad \alpha=\epsilon_i+\epsilon_j\\
0,\qquad & \text{\it if }\quad \alpha\ne\epsilon_i+\epsilon_j
\end{cases}\qquad(\text{mod }\,J),
\end{split}
\end{equation*}
where the last ``$\equiv$" by using the identification w.r.t. modulo
the ideal $J$ as before, and
$a_\ell=\prod\limits_{m=0}^{\ell-1}(\al_j\partial_i-\al_i\pa_j)
(\epsilon_k+m(\al{-}\epsilon_i{-}\epsilon_j)),\
b_\ell=\ell\,(\pa_k-2\pa_{k'})(\al{-}\epsilon_i{-}\epsilon_j)a_{\ell-1}$,
and $a_p=\delta_{ik}-\delta_{jk}$ for
$\alpha=\epsilon_i+\epsilon_j$.

Consequently, by the definition of $d^{(\ell)}$, we get
$d^{(\ell)}((x^{(\alpha)}D_i)^p)=0$ in
$\mathbf{u}(\mathbf{S}(n;\underline{1}))$ for $2\le \ell\le p-1$ and
$0\le\alpha\le\tau$.
\end{proof}

Based on Theorem 3.1, Definition 3.3 and Lemma 3.4, we arrive at
\begin{theorem} Fix two distinguished elements
$h:=D_{kk'}(x^{(\epsilon_k+\epsilon_{k'})})$,
$e:=2D_{kk'}(x^{(2\epsilon_k+\epsilon_{k'})})$ $(1\leq k\neq k'\leq
n)$, there is a noncommutative and noncocummtative Hopf algebra
$(\mathbf{u}_{t,q}(\mathbf{S}(n;\underline{1})),m,\iota,\Delta,S,\varepsilon)$
over $\mathcal{K}[t]_p^{(q)}$ with its algebra structure undeformed,
whose coalgebra structure is given by
\begin{gather*}
\Delta(D_{ij}(x^{(\alpha)}))=D_{ij}(x^{(\alpha)}){\otimes}
(1{-}et)^{\alpha_k{-}\delta_{ik}{-}\delta_{jk}-\a_{k'}{+}\delta_{ik'}{+}\delta_{jk'}}
\tag{36}
\\
\qquad \qquad\qquad
+\sum\limits_{\ell=0}^{p{-}1}{({-}1)}^\ell h^{\lg
\ell\rg}{\otimes}(1{-}et)^{{-}\ell}d^{(\ell)}(D_{ij}(x^{(\alpha)}))t^\ell,
\\
S(D_{ij}(x^{(\alpha)})){=}{-}(1{-}et)^{-\alpha_k{+}\delta_{ik}{+}\delta_{jk}+\a_{k'}{-}\delta_{ik'}{-}\delta_{jk'}}
\cdot\Bigl(\sum\limits_{\ell=0}^{p{-}1}d^{(\ell)}(D_{ij}(x^{(\alpha)}))
\cdot h_1^{\lg \ell\rg}t^\ell\Bigr), \tag{37}\\
\varepsilon(D_{ij}(x^{(\alpha)}))=0, \tag{38}
\end{gather*}
for $0\le\alpha\le\tau$, which is finite dimensional with
$\dim_{\mathcal{K}}\mathbf{u}_{t,q}(\mathbf{S}(n;\underline{1}))=p^{1{+}(n-1)(p^n-1)}$.
\end{theorem}
\begin{proof}
Set $U_{t,q}(\mathbf S(n;\underline 1)):=U(\mathbf S(n;\underline
1))\otimes_{\mathcal K}\mathcal K[t]_p^{(q)}$. Note that the result
of Theorem 3.1, via the base change with $\mathcal K[t]$ instead of
$\mathcal K[t]_p^{(q)}$, is still valid over $U_{t,q}(\mathbf
S(n;\underline 1))$.
 Denote by $I_{t,q}$ the ideal of $U_{t,q}(\mathbf S(n;\underline
1))$ over the ring $\mathcal K[t]_p^{(q)}$ generated by the same
generators of the ideal $I$ in $U(\mathbf S(n;\underline 1))$ via
the base change with $\mathcal K$ replaced by $\mathcal
K[t]_p^{(q)}$. We shall show that $I_{t,q}$ is a Hopf ideal of
$U_{t,q}(\mathbf S(n;\underline 1))$. It suffices to verify that
$\Delta$ and $S$ preserve the generators in $I_{t,q}$ of
$U_{t,q}(\mathbf S(n;\underline 1))$.

\smallskip
(I) \ By Lemmas 2.5, 3.2 \& 3.4 (iii), we obtain
\begin{equation*}
\begin{split}
\Delta((D_{ij}(x^{(\a)}&))^p) =(D_{ij}(x^{(\a)}))^p\otimes
(1{-}et)^{p\,(\alpha_k{-}\a_{k'})}\\
&\qquad\qquad+\sum\limits_{\ell=0}^{\infty} ({-}1)^\ell h^{\lg
\ell\rg}\otimes
(1{-}et)^{{-}\ell}d^{(\ell)}((D_{ij}(x^{(\a)}))^p)t^\ell\
\end{split}\tag{39}
\end{equation*}
\begin{equation*}
\begin{split}
&\equiv(D_{ij}(x^{(\a)}))^p{\otimes}1+\sum\limits_{\ell=0}^{p{-}1}
({-}1)^\ell h^{\lg \ell\rg}{\otimes}
(1{-}et)^{{-}\ell}d^{(\ell)}((D_{ij}(x^{(\a)}))^p)t^\ell\quad
(\text{\rm mod }\, p)
\\
&=(D_{ij}(x^{(\a)}))^p{\otimes}1+1{\otimes}(D_{ij}(x^{(\a)}))^p
+h{\otimes}(1{-}et)^{-1}(\delta_{ik}{-}\delta_{jk})\delta_{\a,\epsilon_i{+}\epsilon_j}et.
\end{split}
\end{equation*}

Hence, when $\alpha\ne\epsilon_i+\epsilon_j$, we get
\begin{equation*}
\begin{split}
\Delta((D_{ij}(x^{(\a)}))^p)&=(D_{ij}(x^{(\a)}))^p\otimes 1+1\otimes
(D_{ij}(x^{(\a)}))^p\\
& \in I_{t,q}\otimes
U_{t,q}(\mathbf{S}(n;\underline{1}))+U_{t,q}(\mathbf{S}(n;\underline{1}))\otimes
I_{t,q};
\end{split}
\end{equation*}
and when $\al=\epsilon_i+\epsilon_j$, by Lemma 3.4 (ii), (28)
becomes
$$
\Delta(D_{ij}(x^{(\epsilon_i+\epsilon_j)}))=D_{ij}(x^{(\epsilon_i+\epsilon_j)}){\otimes}
1+ 1{\otimes} D_{ij}(x^{(\epsilon_i+\epsilon_j)})+h{\otimes}
(1{-}et)^{-1}(\delta_{ik}{-}\delta_{jk})et.
$$
Combining with (39), we obtain
\begin{equation*}
\begin{split}
\Delta((D_{ij}(x^{(\epsilon_i+\epsilon_j)}))^p-D_{ij}(x^{(\epsilon_i+\epsilon_j)})
)&\equiv((D_{ij}(x^{(\epsilon_i+\epsilon_j)}))^p-D_{ij}(x^{(\epsilon_i+\epsilon_j)})
)\otimes 1\\ &\quad +1\otimes
((D_{ij}(x^{(\epsilon_i+\epsilon_j)}))^p-D_{ij}(x^{(\epsilon_i+\epsilon_j)})
)\\
&\in I_{t,q}\otimes
U_{t,q}(\mathbf{S}(n;\underline{1}))+U_{t,q}(\mathbf{S}(n;\underline{1}))\otimes
I_{t,q}.
\end{split}
\end{equation*}

Thereby, we prove that the ideal $I_{t,q}$ is also a coideal of the
Hopf algebra $U_{t,q}(\mathbf{S}(n;\underline{1}))$.

\smallskip
(II) \ By Lemmas 2.5, 3.2 \& 3.4 (iii), we have
\begin{equation*}
\begin{split}
S((D_{ij}(x^{(\a)}))^p) &=-(1{-}et)^{-p(\alpha_k-\a_{k'})}
\sum\limits_{\ell=0}^{\infty} d^{(\ell)}((D_{ij}(x^{(\a)}))^p)\cdot
h_1^{\lg
\ell\rg}t^\ell\\
&\equiv -(D_{ij}(x^{(\a)}))^p-\sum\limits_{\ell=1}^{p-1}
d^{(\ell)}((D_{ij}(x^{(\a)}))^p)\cdot h_1^{\lg \ell\rg}t^\ell \quad
(\text{mod }\,p)\\
&=-(D_{ij}(x^{(\a)}))^p+(\delta_{ik}-\delta_{jk})\delta_{\a,\epsilon_i+\epsilon_j}e\cdot
h_1^{\lg 1\rg} t.
\end{split}\tag{40}
\end{equation*}
Hence, when $\alpha\ne\epsilon_i+\epsilon_j$, we get
$$
S\bigl((D_{ij}(x^{(\a)}))^p\bigr)=-(D_{ij}(x^{(\a)}))^p\in
I_{t,q}.
$$
When $\al=\epsilon_i+\epsilon_j$, by Lemma 3.4 (ii), (29) reads as
$$
S(D_{ij}(x^{(\epsilon_i+\epsilon_j)}))
=-D_{ij}(x^{(\epsilon_i+\epsilon_j)})+(\delta_{ik}-\delta_{jk})
e\cdot h_1^{\lg 1\rg} t.
$$
Combining with (40), we obtain
$$
S\bigl((D_{ij}(x^{(\epsilon_i+\epsilon_j)}))^p-D_{ij}(x^{(\epsilon_i+\epsilon_j)})\bigr)
=-\bigl((D_{ij}(x^{(\epsilon_i+\epsilon_j)}))^p-D_{ij}(x^{(\epsilon_i+\epsilon_j)})\bigr)
\in I_{t,q}.
$$

Thereby, the ideal $I_{t,q}$ is indeed preserved by the antipode $S$
of the quantization $U_{t,q}(\mathbf{S}(n;\underline{1}))$, the same
as in Theorem 3.1.

\smallskip
(III) It is obvious to notice that
$\varepsilon((D_{ij}(x^{(\a)}))^p)=0$ for all $\alpha$ with
$0\le\alpha\le\tau$.

\smallskip
In other words, we prove that $I_{t,q}$ is a Hopf ideal in
 $U_{t,q}(\mathbf{S}(n;\underline{1}))$. We thus obtain the required
 $t$-deformation on
 $\mathbf{u}_{t,q}(\mathbf{S}(n;\underline{1}))$, for the Cartan type
 simple modular restricted Lie algebra of $\mathbf{S}$ type
 --- the special algebra
$\mathbf{S}(n;\underline{1})$.
\end{proof}

\begin{remark} (i) \
Set $f=(1-et)^{-1}$. By Lemma 3.4 \& Theorem 3.5, one gets
$$
[h,f]=f^2-f,\quad h^p=h, \quad f^p=1, \quad
\Delta(h)=h\otimes f+1\otimes h,
$$
where $f$ is a group-like element, and $S(h)=-hf^{-1}$,
$\varepsilon(h)=0$. So the subalgebra generated by $h$ and $f$ is a
Hopf subalgebra of $\mathbf{u}_{t,q}(\mathbf{S}(n;\underline{1}))$,
which is isomorphic to the well-known Radford Hopf algebra over
$\mathcal{K}$ in char $p$ (see \cite{R}).

\smallskip
(ii) \ According to our argument, given a parameter
$q\in\mathcal{K}$, one can specialize $t$ to any root of the
$p$-polynomial $t^p-qt\in\mathcal{K}[t]$ in a split field of
$\mathcal{K}$. For instance, if take $q=1$, then one can specialize
$t$ to any scalar in $\mathbb{Z}_p$. If set $t=0$, then we get the
original standard Hopf algebra structure of
$\mathbf{u}(\mathbf{S}(n;\underline{1}))$. In this way, we indeed
get a new Hopf algebra structure over the same restricted universal
enveloping algebra $\mathbf{u}(\mathbf{S}(n;\underline{1}))$ over
$\mathcal{K}$ under the assumption that $\mathcal{K}$ is
algebraically closed, which has the new coalgebra structure induced
by Theorem 3.5, but has dimension $p^{(n-1)(p^n-1)}$.
\end{remark}

\section{More quantizations}
In this section, we can get more Drinfel'd twists by considering the
products of some pairwise different {\it basic Drinfel'd twists} as
stated in Remark 2.8. By the same argument as in Theorem 2.4, one
can get many more new complicated quantizations not only over the
$U(\mathbf{S}_{\mathbb{Z}}^+)[[t]]$, but over the
$\mathbf{u}_{t,q}(\mathbf{S}(n;\underline{1}))$ as well. Moreover,
we prove that the twisted structures given by some products of
pairwise different {\it basic Drinfel'd twists} with different
length are nonisomorphic.

\subsection{More Drinfel'd twists}
We consider the products of pairwise different and mutually
commutative basic Drinfel'd twists. Note that $[\mathcal{F}(i,j),
\mathcal{F}(k,m)]=0$ for $i\neq k,m$ and $j\neq k$. This fact,
according to the definition of $\mathcal{F}(k,m)$, implies the
commutative relations in the case $i\neq k, m$ and $j\neq k$:
\begin{equation*}
\begin{split}
(\mathcal{F}(k,m)\otimes 1)(\Delta_0\otimes\text{\rm Id})
(\mathcal{F}(i,j))&=(\Delta_0\otimes\text{\rm Id})
(\mathcal{F}(i,j))(\mathcal{F}(k,m)\otimes 1),\\
(1\otimes \mathcal{F}(k,m))(\text{\rm Id}\otimes\Delta_0)
(\mathcal{F}(i,j))&=(\text{\rm Id}\otimes\Delta_0)
(\mathcal{F}(i,j))(1\otimes\mathcal{F}(k,m)),
\end{split}\tag{41}
\end{equation*}
which give rise to the following property.
\begin{theorem}
$\mathcal{F}(i,j)\mathcal{F}(k,m)(i \neq k,m; j\neq k)$ is still a
Drinfel'd twist on $U(\mathbf{S}_{\mathbb{Z}}^+)[[t]]$.
\end{theorem}
\begin{proof}
Note that $\Delta_0\otimes\text{\rm id}$, $\text{\rm
id}\otimes\Delta_0$, $\varepsilon_0\otimes\text{\rm id}$ and
$\text{\rm id}\otimes\varepsilon_0$ are algebraic homomorphisms.
According to Lemma 1.4, it suffices to check that
\begin{equation*}
\begin{split}
(\mathcal{F}(i,j)\mathcal{F}(k,m)\otimes 1)&(\Delta_0\otimes
\text{\rm Id})(\mathcal{F}(i,j)\mathcal{F}(k,m))\\
&= (1\otimes \mathcal{F}(i,j)\mathcal{F}(k,m))(\text{\rm
Id}\otimes\Delta_0) (\mathcal{F}(i,j)\mathcal{F}(k,m)).
\end{split}
\end{equation*}

Using $(41)$, we have
\begin{equation*}
\begin{split}
\text{LHS}&=(\mathcal{F}(i,j)\otimes 1)(\mathcal{F}(k,m)\otimes
1)(\Delta_0\otimes \text{\rm Id})(\mathcal{F}(i,j))(\Delta_0\otimes
\text{\rm
Id})(\mathcal{F}(k,m))\\
&=(\mathcal{F}(i,j)\otimes 1)(\Delta_0\otimes \text{\rm
Id})(\mathcal{F}(i,j))(\mathcal{F}(k,m)\otimes 1)(\Delta_0\otimes
\text{\rm
Id})(\mathcal{F}(k,m))\\
&=(1\otimes\mathcal{F}(i,j) )(\text{\rm Id}\otimes\Delta_0
)(\mathcal{F}(i,j))(1\otimes \mathcal{F}(k,m))(\text{\rm
Id}\otimes\Delta_0)(\mathcal{F}(k,m))\\
&=(1\otimes\mathcal{F}(i,j) )(1\otimes \mathcal{F}(k,m))(\text{\rm
Id}\otimes\Delta_0 )(\mathcal{F}(i,j))(\text{\rm
Id}\otimes\Delta_0)(\mathcal{F}(k,m))=\text{RHS}.
\end{split}
\end{equation*}

This completes the proof.
\end{proof}
More generally, we have the following
 \begin{coro}
Let $\mathcal{F}(i_1,j_1),\cdots, \mathcal{F}(i_m,j_m)$ be $m$
pairwise different basic Drinfel'd twists and
 $[\mathcal{F}(i_k,j_k), \mathcal{F}(i_s,j_s)]=0$ for all $1\leq k\neq s\leq
 m$. Then $\mathcal{F}(i_1,j_1) \cdots \mathcal{F}(i_m,j_m)$
 is still
a Drinfel'd twist.
\end{coro}

We denote $\mathcal{F}_m=\mathcal{F}(i_1,j_1)\cdots
\mathcal{F}(i_m,j_m)$ and the length of  $\mathcal{F}(i_1,j_1)\cdots
\mathcal{F}(i_m,j_m)$ as $m$. These twists lead to more
quantizations.

\subsection{More quantizations}
We consider the modular reduction process for the quantizations of
$U(\mathbf{S}^+)[[t]]$ arising from those products of some pairwise
different and mutually commutative basic Drinfel'd twists. We will
then get lots of new families of noncommutative and noncocommutative
Hopf algebras of dimension $p^{1{+}(n-1)(p^n-1)}$ with indeterminate
$t$ or of dimension $p^{(n-1)(p^n-1)}$ with specializing $t$ into a
scalar in $\mathcal{K}$.

Let $A(k,k')_\ell$,\, $B(k,k')_\ell$ and $A(m,m')_n$,\, $B(m,m')_n$
denote the coefficients of the corresponding quantizations of
$U(\mathbf{S}^+_{\mathbb{Z}})$ over
$U(\mathbf{S}^+_{\mathbb{Z}})[[t]]$ given by Drinfel'd twists
$\mathcal{F}(k,k')$ and $\mathcal{F}(m,m')$ respectively as in
Corollary 2.7. Note that $A(k,k')_0=A(m,m')_0$ $=1$,
$A(k,k')_{-1}=A(m,m')_{-1}=0$.

Set
\begin{equation*}
\begin{split}
\pa(m, m';k, k')_{\ell,n}&:=A(m,m')_nA(k,k')_\ell\pa-A(m,m')_nB(k,k')_\ell(\pa_k{-}2\pa_{k'})\\
&\qquad\qquad\quad\,-A(k,k')_\ell B(m,m')_n(\pa_{m}{-}2\pa_{m'}).
\end{split}
\end{equation*}

\begin{lemm}
Fix distinguished elements $h(k,k')=\pa_k{-}\pa_{k'}$,
$e(k,k')=x^{\epsilon_k}(\pa_k-2\pa_{k'})$ $(1\le k\neq k'\le n)$
 and $h(m,m')=\pa_{m}{-}\pa_{m'}$,
$e(m,m')=x^{\epsilon_m}(\pa_{m}{-}2\pa_{m'})$ $(1\le m\neq m'\le n)$
with $k\neq m,m'$ and $k'\neq m$, the corresponding quantization of
$U(\mathbf{S}^+_{\mathbb{Z}})$ over
$U(\mathbf{S}^+_{\mathbb{Z}})[[t]]$ by Drinfel'd twist
$\mathcal{F}=\mathcal{F}(m,m')\mathcal{F}(k,k')$ with the product
undeformed is given by
\begin{gather*}
\Delta(x^{\alpha}\pa )= x^{\alpha}\pa \otimes
\bigl(1{-}e(k,k')t\bigr)^{\alpha_k{-}\a_{k'}}\bigl(1{-}e(m,m')t\bigr)^{\alpha_m{-}\a_{m'}}
\tag{42} \\\qquad\qquad\qquad\
+\sum\limits_{n,\ell=0}^{\infty}{({-}1)}^{n+\ell} h(k,k')^{\lg
\ell\rg}\cdot h(m,m')^{\lg n \rg}\otimes \bigl(1{-}e(k,k')t\bigr)^{{-}\ell}\\
\qquad\qquad\qquad\qquad\quad \cdot\,\bigl(1{-}e(m,m')t\bigr)^{{-}n}
x^{\alpha+\ell\epsilon_k+n\epsilon_m} \pa(m,m';k,k')_{\ell,n}
t^{n+\ell},\\
S(x^{\alpha}\pa)={-}\bigl(1{-}e(k,k')t\bigr)^{-\alpha_k{+}\a_{k'}}
\bigl(1{-}e(m,m')t\bigr)^{-\alpha_m{+}\a_{m'}}\cdot\tag {43}\\
\qquad\qquad\qquad\qquad\ \cdot \sum\limits_{n,\ell=0}^{\infty}
x^{\alpha+\ell\epsilon_k+n\epsilon_m} \pa(m,m';k,k')_{\ell,n} \cdot
h(m,m')_{1}^{\lg n\rg}h(k,k')_{1}^{\lg \ell\rg} t^{n+\ell},\\
 \varepsilon(x^{\alpha}\pa
)=0,  \tag{44}
\end{gather*}
for  $x^\alpha\pa \in \mathbf{S}^+_{\mathbb{Z}}$.
\end{lemm}
\begin{proof} Using Corollary 2.7, we get
\begin{equation*}
\begin{split}
\Delta(x^{\alpha}\pa)&=\mathcal{F}(m,m')\mathcal{F}(k,k')
\Delta_0(x^{\alpha}\pa)\mathcal{F}(k,k')^{-1}\mathcal{F}(m,m')^{-1}\\
&=\mathcal{F}(m,m')\Bigl( x^{\alpha}\pa\otimes
\bigl(1{-}e(k,k')t\bigr)^{\alpha_k{-}\a_{k'}}\\
&\qquad+\sum\limits_{\ell=0}^{\infty}{({-}1)}^\ell h(k,k')^{\lg
\ell\rg}\otimes\bigl(1{-}e(k,k')t\bigr)^{{-}\ell}\cdot
x^{\alpha{+}\ell\epsilon_k}\bigl(A(k,k')_\ell \pa \\
&\qquad\qquad -B(k,k')_\ell(\pa_k{-}2\pa_{k'})\bigr)t^\ell\Bigr)
\mathcal{F}(m,m')^{-1}.
\end{split}
\end{equation*}
Using (19) and Lemma 2.1, we get
\begin{equation*}
\begin{split}
\mathcal{F}(m,m')&\Bigl( x^{\alpha}\pa\otimes
\bigl(1{-}e(k,k')t\bigr)^{\alpha_k{-}\a_{k'}}\Bigr) \mathcal{F}(m,m')^{-1}\\
&=\mathcal{F}(m,m')\Bigl( x^{\alpha}\pa\otimes 1 \Bigr)
\mathcal{F}(m,m')^{-1}
\Bigl(1\otimes\bigl(1{-}e(k,k')t\bigr)^{\alpha_k{-}\a_{k'}}\Bigr)
\\
&=\mathcal{F}(m,m')\mathcal{F}(m,m')^{-1}_{\a_{m'}{-}\alpha_m}
\Bigl( x^{\alpha}\pa \otimes 1 \Bigr)\Bigl(1\otimes
\bigl(1{-}e(k,k')t\bigr)^{\alpha_k{-}\a_{k'}}\Bigr)\\
  &=\Bigl(1\otimes\bigl(1{-}e(m,m')t\bigr)^{\alpha_m{-}\a_{m'}}\Bigr)
 \Bigl( x^{\alpha}\pa \otimes 1 \Bigr)
 \Bigl(1\otimes\bigl(1{-}e(k,k')t\bigr)^{\alpha_k{-}\a_{k'}}\Bigr)\\
&=x^{\alpha}\pa \otimes
\bigl(1{-}e(k,k')t\bigr)^{\alpha_k{-}\a_{k'}}\bigl(1{-}e(m,m')t\bigr)^{\alpha_m{-}\a_{m'}}.
\end{split}
\end{equation*}
Using (21), we  have
\begin{equation*}
\begin{split}
\mathcal{F}&(m,m')\Bigl(\sum\limits_{\ell=0}^{\infty}{({-}1)}^\ell
h(k,k')^{\lg \ell\rg}\otimes \bigl(1{-}e(k,k')t\bigr)^{{-}\ell} \\
&\qquad \cdot x^{\alpha{+}\ell\epsilon_k}\bigl(A(k,k')_\ell
\pa-B(k,k')_\ell(\pa_k{-}2\pa_{k'})\bigr)t^\ell\Bigr) \mathcal{F}(m,m')^{-1}\\
&=\sum\limits_{\ell=0}^{\infty}{({-}1)}^\ell h(k,k')^{\lg
\ell\rg}\otimes \bigl(1{-}e(k,k')t\bigr)^{{-}\ell}\\
 &\qquad \cdot
\mathcal{F}(m,m')\Bigl(1 \otimes
x^{\alpha{+}\ell\epsilon_k}\bigl(A(k,k')_\ell
\pa-B(k,k')_\ell(\pa_k{-}2\pa_{k'})\bigr)\Bigr)t^\ell \mathcal{F}(m,m')^{-1}\\
&=\sum\limits_{n,\ell=0}^{\infty}{({-}1)}^\ell h(k,k')^{\lg
\ell\rg}\otimes \bigl(1{-}e(k,k')t\bigr)^{{-}\ell}\cdot
\mathcal{F}(m,m')F(m,m')_{n} \Bigl( h(m,m')^{\lg n \rg}\otimes\\
&\qquad\quad
 x^{\alpha+\ell\epsilon_k+n\epsilon_m}\pa(m,m';k,k')_{\ell,n}
t^{n+\ell}\Bigr)
\\
&=\sum\limits_{n,\ell=0}^{\infty}{({-}1)}^\ell h(k,k')^{\lg
\ell\rg}h(m,m')^{\lg n \rg}\otimes
\bigl(1{-}e(k,k')t\bigr)^{{-}\ell}\bigl(1{-}e(m,m')t\bigr)^{{-}n}\\
& \qquad \cdot
x^{\alpha+\ell\epsilon_k+n\epsilon_m}\pa(m,m';k,k')_{\ell,n}t^{n+\ell}.
\end{split}
\end{equation*}

For $k\neq m, m'$ and $k'\neq m$, by the definitions of $v$ and $u$,
we get
\begin{equation*}
\begin{split}
v&=v(k,k')v(m,m')=v(m,m')v(k,k'),\\
u&=u(m,m')u(k,k')=u(k,k')u(m,m').
\end{split}
\end{equation*}

Note $ u(m,m')h(k,k')=h(k,k')u(m,m')$,
$v(m,m')e(k,k')=e(k,k')v(m,m')$. By Corollary 2.7 and (20), we have
\begin{equation*}
\begin{split}
S(x^{\alpha}\pa)&=-v\cdot x^{\alpha}\pa\cdot u \\
&=-v(m,m')v(k,k')\cdot x^{\alpha}\pa\cdot  u(k,k')u(m,m')\\
&=v(m,m')\cdot\Bigl({-}\bigl(1{-}e(k,k')t\bigr)^{-\alpha_k{+}\a_{k'}}\cdot
\Bigl(\sum\limits_{\ell=0}^{\infty}
x^{\alpha{+}\ell\epsilon_k}\bigl(A(k,k')_\ell\pa \\
&\qquad -B(k,k')_\ell(\pa_k{-}2\pa_{k'})\bigr)\cdot h(k,k')_1^{\lg
\ell\rg}t^\ell\Bigr)\Bigr)\cdot u(m,m')
\\
&={-}\bigl(1{-}e(k,k')t\bigr)^{-\alpha_k{+}\a_{k'}}\cdot
v(m,m')u(m,m')_{\alpha_m-\a_{m'}} \\
&\qquad \cdot \sum\limits_{n,\ell=0}^{\infty}
x^{\alpha+\ell\epsilon_k+n\epsilon_m}\pa(m,m';k,k')_{\ell,n}
 \cdot h(m,m')_{1}^{\lg
n\rg}h(k,k')_{1}^{\lg
\ell\rg} t^{n+\ell}\\
\end{split}
\end{equation*}
\begin{equation*}
\begin{split}
&={-}\bigl(1{-}e(k,k')t\bigr)^{-\alpha_k{+}\a_{k'}}\bigl(1{-}e(m,m')t\bigr)^{-\alpha_m{+}\a_{m'}}
\\
 &\qquad \cdot \sum\limits_{n,\ell=0}^{\infty}
x^{\alpha+\ell\epsilon_k+n\epsilon_m}\pa(m,m';k,k')_{\ell,n}
 \cdot h(m,m')_{1}^{\lg
n\rg}h(k,k')_{1}^{\lg \ell\rg} t^{n+\ell}.
\end{split}
\end{equation*}

This completes the proof.
\end{proof}

Set $\a
(k,k')=\a_k{-}\delta_{ik}{-}\delta_{jk}{-}\a_{k'}{+}\delta_{ik'}{+}\delta_{jk'}$
and $d_{kk'}^{(\ell)}=\frac{1}{\ell!}(\text{\rm
ad}\,e(k,k'))^{\ell}$. Write coefficients $\bar A_\ell$, $\bar
B_\ell$, $ A_\ell$ in Theorem 3.1 as $\bar A(k,k')_\ell$, $\bar
B(k,k')_\ell$, $ A(k,k')_\ell$, respectively. Set
\begin{equation*}
\begin{split}
D_{ij}(m, m'; &k, k')_{\ell, n}:=\bar{A}(k,k')_\ell \bar{A}(m,m')_n
D_{ij}(x^{(\alpha{+}\ell\epsilon_k+n\epsilon_m)})\\
&\,+\bar{B}(k,k')_\ell \bar{A}(m,m')_n (\delta_{ik}D_{k'j}
+\delta_{jk}D_{ik'})(x^{(\alpha{+}(\ell-1)\epsilon_k+n\epsilon_m+\epsilon_{k'})})
\\
&\,+\bar{A}(k,k')_\ell \bar{B}(m,m')_n
(\delta_{im}D_{k'j}+\delta_{jm}D_{ik'})
(x^{(\alpha{+}\ell\epsilon_k+(n-1)\epsilon_m+\epsilon_{k'})}).
\end{split}
\end{equation*}

Using Lemma 4.3, we get a new quantization of
$U(\mathbf{S}(n;\underline{1}))$ over
$U_t(\mathbf{S}(n;\underline{1}))$ by Drinfel'd twist
$\mathcal{F}=\mathcal{F}(m,m')\mathcal{F}(k,k')$ as follows.

\begin{lemm}
Fix distinguished elements
$h(k,k')=D_{kk'}(x^{(\epsilon_k+\epsilon_{k'})})$,
$e(k,k')=2D_{kk'}(x^{(2\epsilon_k+\epsilon_{k'})})$;
$h(m,m')=D_{mm'}(x^{(\epsilon_m+\epsilon_{m'})})$,
$e(m,m')=2D_{mm'}(x^{(2\epsilon_m+\epsilon_{m'})})$ with $k\neq
m,m'; k'\neq m$, the corresponding quantization of
$U(\mathbf{S}(n;\underline{1}))$ on
$U_t(\mathbf{S}(n;\underline{1}))$ $($also on
$U(\mathbf{S}(n;\underline{1}))[[t]])$ with the product undeformed
is given by
\begin{gather*}
\Delta(D_{ij}(x^{(\alpha)}))=D_{ij}(x^{(\alpha)})\otimes
\bigl(1{-}e(k,k')t\bigr)^{\alpha(k,k')}\bigl(1{-}e(m,m')t\bigr)^{\alpha(m,m')}\tag{45}\\
\qquad\qquad\qquad\qquad\,+\sum\limits_{n,\ell=0}^{p{-}1}{({-}1)}^{n+\ell}h(k,k')^{\lg
\ell\rg}h(m,m')^{\lg
n\rg}\otimes \bigl(1{-}e(k,k')t\bigr)^{{-}\ell}\cdot\\
\qquad\qquad\qquad\qquad\qquad\qquad
\cdot\bigl(1{-}e(m,m')t\bigr)^{{-}n}D_{ij}(m, m'; k, k')_{\ell,
n}t^{n+\ell},
\\
S(D_{ij}(x^{(\alpha)}))={-}\bigl(1{-}e(k,k')t\bigr)^{-\alpha(k,k')}
\bigl(1{-}e(m,m')t\bigr)^{-\alpha(m,m')}\cdot \tag{46} \\
\qquad\qquad\qquad\qquad\qquad
\cdot\,\Bigl(\sum\limits_{n,\ell=0}^{p{-}1}D_{ij}(m,m';k,k')_{\ell,n}
h(k,k')_1^{\lg \ell\rg}h(m,m')_1^{\lg n\rg}t^{n+\ell}\Bigr),
\\
\varepsilon(D_{ij}(x^{(\alpha)}))=0, \tag{47}
\end{gather*}
where  $0\le \alpha \le \tau$.
\end{lemm}

For the further discussion, we need two lemmas below about the
quantization of $U(\mathbf{S}(n;\underline{1}))$ over
$U(\mathbf{S}(n;\underline{1}))[[t]]$ in Lemma 4.4.
\begin{lemm} For $s\ge 1$, one has
\begin{gather*}
\Delta((D_{ij}(x^{(\alpha)}))^s)=\sum_{0\le j\le s\atop n, \ell\ge
0}\dbinom{s}{j}({-}1)^{n+\ell}(D_{ij}(x^{(\alpha)}))^jh(k,k')^{\langle
\ell\rangle}h(m,m')^{\langle
n\rangle}\otimes\tag{\text{\rm i}} \\
\qquad\qquad\qquad\qquad\quad
\bigl(1{-}e(k,k')t\bigr)^{j\al(k,k'){-}\ell}\bigl(1{-}e(m,m')t\bigr)^{j\al(m,m'){-}n}\cdot
\\
\qquad\qquad\qquad\qquad\qquad \cdot\,
d_{mm'}^{(n)}d_{kk'}^{(\ell)}((D_{ij}(x^{(\alpha)}))^{s{-}j})t^{\ell+n}.
\\
S((D_{ij}(x^{(\alpha)}))^s)=(-1)^s\bigl(1{-}e(m,m')t\bigr)^{-s\al(m,m')}
\bigl(1{-}e(k,k')t\bigr)^{-s\al(k,k')}\cdot\tag{\text{\rm ii}} \\
\qquad\qquad\qquad\qquad\quad\; \cdot
\Bigl(\sum\limits_{n,\ell=0}^{\infty}
d_{mm'}^{(n)}d_{kk'}^{(\ell)}((D_{ij}(x^{(\alpha)}))^s)
h(k,k')_1^{\lg \ell\rg} h(m,m')_1^{\lg n \rg}t^{n+\ell}\Bigr).
\end{gather*}
\end{lemm}
\begin{proof} By Lemma 2.5, (21), (19) and Lemma 1.6, we obtain
\begin{equation*}
\begin{split}
&\Delta((D_{ij}(x^{(\alpha)}))^s)
=\mathcal{F}\Bigl(D_{ij}(x^{(\alpha)})\otimes
1+1\otimes D_{ij}(x^{(\alpha)})\Bigr)^s\mathcal{F}^{-1}\\
&=\mathcal{F}(m,m')\Bigl(\sum_{0\le j\le s\atop
\ell\ge0}\dbinom{s}{j}({-}1)^\ell(D_{ij}(x^{(\alpha)}))^jh(k,k')^{\langle
\ell\rangle}\otimes\bigl(1{-}e(k,k')t\bigr)^{j\al(k,k'){-}\ell}\\
& \quad
\cdot d_{kk'}^{(\ell)}((D_{ij}(x^{(\alpha)}))^{s{-}j})t^\ell\Bigr)\mathcal{F}(m,m')^{-1}\\
&=\mathcal{F}(m,m')\Bigl(\sum_{0\le j\le s\atop
\ell\ge0}\dbinom{s}{j}({-}1)^\ell\bigl((D_{ij}(x^{(\alpha)}))^j{\otimes}
1\bigr) \bigl(h(k,k')^{\langle
\ell\rangle}{\otimes}\bigl(1{-}e(k,k')t\bigr)^{j\al(k,k'){-}\ell}\bigr)\\
& \quad
\cdot \bigl(1\otimes d_{kk'}^{(\ell)}((D_{ij}(x^{(\alpha)}))^{s{-}j})t^\ell\bigr)\Bigr)\mathcal{F}(m,m')^{-1}\\
&=\mathcal{F}(m,m')\sum_{0\le j\le s\atop n, \ell\ge
0}\dbinom{s}{j}({-}1)^{n+\ell}\bigl((D_{ij}(x^{(\alpha)}))^j\otimes
1\bigr)\mathcal{F}(m,m')_n^{-1}h(k,k')^{\langle
\ell\rangle}h(m,m')^{\langle n\rangle}\\  & \quad
\otimes\bigl(1{-}e(k,k')t\bigr)^{j\al(k,k'){-}\ell}
d_{mm'}^{(n)}d_{kk'}^{(\ell)}((D_{ij}(x^{(\alpha)}))^{s{-}j})t^{\ell+n}\\
&=\sum_{0\le j\le s\atop n, \ell\ge
0}\dbinom{s}{j}({-}1)^{n+\ell}\mathcal{F}(m,m')\mathcal{F}(m,m')_{n-j\al(m,m')}^{-1}
\bigl((D_{ij}(x^{(\alpha)}))^j\otimes 1\bigr)h(k,k')^{\langle \ell\rangle}\\
& \quad \cdot h(m,m')^{\langle n\rangle}
\otimes\bigl(1{-}e(k,k')t\bigr)^{j\al(k,k'){-}\ell}
d_{mm'}^{(n)}d_{kk'}^{(\ell)}((D_{ij}(x^{(\alpha)}))^{s{-}j})t^{\ell+n}\\
&=\sum_{0\le j\le s\atop n, \ell\ge
0}\dbinom{s}{j}({-}1)^{n+\ell}(D_{ij}(x^{(\alpha)}))^jh(k,k')^{\langle
\ell\rangle}h(m,m')^{\langle
n\rangle}\otimes\bigl(1{-}e(k,k')t\bigr)^{j\al(k,k'){-}\ell}
\\  & \quad \cdot\bigl(1{-}e(m,m')t\bigr)^{j\al(m,m'){-}n}
d_{mm'}^{(n)}d_{kk'}^{(\ell)}((D_{ij}(x^{(\alpha)}))^{s{-}j})t^{\ell+n}.
\end{split}
\end{equation*}
Again by (20) and Lemma 1.6,
\begin{equation*}
\begin{split}
S((D_{ij}(x^{(\alpha)}))^s)&=u^{-1}S_0((D_{ij}(x^{(\alpha)}))^s)\,u\\
&=(-1)^s v\cdot (D_{ij}(x^{(\alpha)}))^s\cdot u \\ &=(-1)^s
v(m,m')\Bigl(\bigl(1{-}e(k,k')t\bigr)^{-s\al(k,k')}
\\ & \quad \cdot\Bigl(\sum\limits_{\ell=0}^{\infty}
d_{kk'}^{(\ell)}((D_{ij}(x^{(\alpha)}))^s)\cdot h(k,k')_1^{\lg
\ell\rg}t^\ell\Bigr)\Bigr)u(m,m')
\\
 &=(-1)^s
v(m,m')u(m,m')_{s\al(m,m')}\bigl(1{-}e(k,k')t\bigr)^{-s\al(k,k')}\\
&\quad \cdot \Bigl(\sum\limits_{n,\ell=0}^{\infty}
d_{mm'}^{(n)}d_{kk'}^{(\ell)}((D_{ij}(x^{(\alpha)}))^s)\cdot
h(k,k')_1^{\lg
\ell\rg} h(m,m')_1^{\lg n \rg}t^{n+\ell}\Bigr) \\
&=(-1)^s
\bigl(1{-}e(m,m')t\bigr)^{-s\al(m,m')}\bigl(1{-}e(k,k')t\bigr)^{-s\al(k,k')}\\
& \quad \cdot \Bigl(\sum\limits_{n,\ell=0}^{\infty}
d_{mm'}^{(n)}d_{kk'}^{(\ell)}((D_{ij}(x^{(\alpha)}))^s)\cdot
h(k,k')_1^{\lg \ell\rg} h(m,m')_1^{\lg n \rg}t^{n+\ell}\Bigr).
\end{split}
\end{equation*}

This completes the proof.
\end{proof}

\begin{lemm} \ Set $e(k,k')=2D_{kk'}(x^{(2\epsilon_k+\epsilon_{k'})})$, $e(m,m')=2D_{mm'}
(x^{(2\epsilon_m+\epsilon_{m'})})$,

\smallskip
$d_{kk'}^{(\ell)}=\frac{1}{\ell!}(\text{\rm ad}\,e(k,k'))^\ell$ and
$d_{mm'}^{(n)}=\frac{1}{n!}(\text{\rm ad}\,e(m,m'))^n$. Then

\smallskip
 $(\text{\rm i})$ \
$d_{mm'}^{(n)}d_{kk'}^{(\ell)}(D_{ij}(x^{(\al)}))=D_{ij}(m,m';k,k')_{\ell,n}$,

\smallskip
\hskip0.6cm where $D_{ij}(m,m';k,k')_{\ell,n}$ as in Lemma 4.4.

\smallskip
$(\text{\rm ii})$ \
$d_{mm'}^{(n)}d_{kk'}^{(\ell)}(D_{ij}(x^{(\epsilon_i+\epsilon_j)}))
=\delta_{\ell,0}\delta_{n,0}D_{ij}(x^{(\epsilon_i+\epsilon_j)})
-\delta_{n,0}\delta_{1,\ell}(\delta_{ik}{-}\delta_{jk})e(k,k')$

\smallskip
\hskip5cm
$-\,\delta_{\ell,0}\delta_{1,n}(\delta_{im}{-}\delta_{jm})e(m,m')$.

\smallskip
$(\text{\rm iii})$ \
$d_{mm'}^{(n)}d_{kk'}^{(\ell)}((D_{ij}(x^{(\al)}))^p)=\delta_{\ell,0}\delta_{n,0}(D_{ij}(x^{(\al)}))^p
-\delta_{n,0}\delta_{1,\ell}
(\delta_{ik}{-}\delta_{jk})\delta_{\alpha,\epsilon_i+\epsilon_j}\cdot$

\smallskip
\hskip5cm $\cdot e(k,k')-\delta_{\ell,0}\delta_{1,n}
(\delta_{im}{-}\delta_{jm})\delta_{\alpha,\epsilon_i+\epsilon_j}e(m,m')$.
\end{lemm}
\begin{proof} (i) For $0\le \alpha\le\tau$, using (17), we
obtain
\begin{equation*}
\begin{split}
d_{mm'}^{(n)}&d_{kk'}^{(\ell)}(D_{ij}(x^{(\al)}))\\
&=d_{mm'}^{(n)}d_{kk'}^{(\ell)}\Bigl(\frac{1}{\a !}
x^{\a-\epsilon_i-\epsilon_j}(\a_j\pa_i-\a_i\pa_j)\Bigr)\\
&=d_{mm'}^{(n)}\Bigl(\frac{1}{\a
!}x^{\a-\epsilon_i-\epsilon_j+\ell\epsilon_k}\bigl(A(k,k')_\ell
(\a_j\pa_i{-}\a_i\pa_j)-B(k,k')_\ell(\pa_k{-}2\pa_{k'})\bigr)\Bigr)\\
&=\frac{1}{\a
!}x^{\a-\epsilon_i-\epsilon_j+\ell\epsilon_k+n\epsilon_m}\bigl(A(k,k')_\ell
A(m,m')_n (\a_j\pa_i{-}\a_i\pa_j)\\
& \quad -A(m,m')_n
B(k,k')_\ell(\pa_k{-}2\pa_{k'})-
A(k,k')_\ell B(m,m')_n(\pa_{k'}{-}2\pa_{k'})\bigr)\\
&=D_{ij}(m,m';k,k')_{\ell,n}.
\end{split}
\end{equation*}

(ii), (iii) may be proved directly using  Lemma 3.4.
\end{proof}
Using Lemmas 3.2, 3.4, 4.5 \&  4.6, we get a new Hopf algebra
structure over the same restricted universal enveloping algebra
$\mathbf{u}(\mathbf{S}(n;\underline{1}))$ over $\mathcal{K}$ by the
products of two different and commutative basic Drinfel'd twists.
\begin{theorem} Fix two distinguished elements
$h(k,k'):$ $=D_{kk'}(x^{(\epsilon_k+\epsilon_{k'})})$,
$e(k,k'):=2D_{kk'}(x^{(2\epsilon_k+\epsilon_{k'})})$ $(1\le k\neq
k'\le n)$ and $h(m,m'):=D_{mm'}(x^{(\epsilon_m+\epsilon_{m'})})$,
$e(m,m'):=2D_{mm'}(x^{(2\epsilon_m+\epsilon_{m'})})$ $(1\le m\neq
m'\le n)$ with $k\neq m, m'; k'\neq m$, there is a noncommutative
and noncocummtative Hopf algebra
$(\mathbf{u}_{t,q}(\mathbf{S}(n;\underline{1})),m,\iota,$
$\Delta,S,\varepsilon)$ over $\mathcal{K}[t]_p^{(q)}$ with the
product undeformed, whose coalgebra structure is given by
\begin{gather*}
\Delta(D_{ij}(x^{(\al)}))=D_{ij}(x^{(\al)})\otimes
\bigl(1{-}e(k,k')t\bigr)^{\alpha(k,k')}\bigl(1{-}e(m,m')t\bigr)^{\alpha(m,m')}  \tag{48}\\
\qquad\qquad\qquad\qquad+
\sum\limits_{n,\ell=0}^{p-1}(-1)^{\ell+n}h(k,k')^{\lg \ell\rg}
h(m,m')^{\lg n \rg}\otimes \bigl(1{-}e(k,k')t\bigr)^{-\ell}\cdot\\
\qquad\qquad\qquad\qquad\qquad\qquad\cdot\, \bigl(1{-}e(m,m')t\bigr)^{-n}d_{kk'}^{(\ell)}d_{mm'}^{(n)}(D_{ij}(x^{(\al)}))t^{\ell+n},\\
S(D_{ij}(x^{(\al)}))=-\bigl(1{-}e(k,k')t\bigr)^{-\alpha(k,k')}
\bigl(1{-}e(m,m')t\bigr)^{-\alpha(m,m')}\cdot \tag{49}\\
\qquad\qquad\qquad\qquad\qquad\cdot\,
\Bigl(\sum\limits_{n,\ell=0}^{p-1}d_{kk'}^{(\ell)}d_{mm'}^{(n)}(D_{ij}(x^{(\al)}))
 h(k,k')_1^{\lg \ell\rg}h(m,m')_1^{\lg n\rg}t^{\ell+n}\Bigr),\\
\varepsilon(D_{ij}(x^{(\al)}))=0,  \tag{50}
\end{gather*}
where $0\le\alpha\le\tau$, and
$\dim_{\mathcal{K}}\mathbf{u}_{t,q}(\mathbf{S}(n;\underline{1}))=p^{1{+}(n-1)(p^n-1)}$.
\end{theorem}
\begin{proof}
Let $I_{t,q}$ denote the ideal of
$(U_{t,q}(\mathbf{S}(n;\underline{1})), m,
\iota,\Delta,S,\varepsilon)$ over the ring $\mathcal K[t]_p^{(q)}$
generated by the same generators as in $I$ ($q\in\mathcal{K}$).
Observe that the result in Lemma 4.4, via the base change with
$\mathcal K[t]$ replaced by $\mathcal K[t]_p^{(q)}$, is still valid
for $U_{t,q}(\mathbf{S}(n;\underline{1}))$.

In what follows, we shall show that $I_{t,q}$ is a Hopf ideal of
$U_{t,q}(\mathbf{S}(n;\underline{1}))$. To this end, it suffices to
verify that $\Delta$ and $S$ preserve the generators of $I_{t,q}$.

\smallskip
(I) \ By Lemmas 4.5, 3.2, 3.4 \& 4.6, we obtain
\begin{equation*}
\begin{split}
\Delta((D_{ij}(&x^{(\al)}))^p)=(D_{ij}(x^{(\al)}))^p\otimes
(1{-}e(k,k')t)^{p\,\alpha(k,k')}(1{-}e(m,m')t)^{p\,\alpha(m,m')}\\
&\quad+\sum\limits_{n,\ell=0}^{\infty} ({-}1)^{n+\ell} h(k,k')^{\lg
\ell\rg} h(m,m')^{\lg n \rg}\otimes
(1{-}e(k,k')t)^{{-}\ell}\\
& \quad \cdot(1{-}e(m,m')t)^{{-}n}
d_{mm'}^{(n)}d_{kk'}^{(\ell)}((D_{ij}(x^{(\al)}))^p)t^{n+\ell}
\\
&\equiv(D_{ij}(x^{(\al)}))^p{\otimes}1{+}\sum\limits_{n,\ell=0}^{p{-}1}
({-}1)^{n+\ell} h(k,k')^{\lg \ell\rg} h(m,m')^{\lg n
\rg}{\otimes}(1{-}e(k,k')t)^{{-}\ell}
\\ & \quad \cdot
(1{-}e(m,m')t)^{{-}n}d_{mm'}^{(n)}d_{kk'}^{(\ell)}((D_{ij}(x^{(\al)}))^p)t^{n+\ell}
\quad (\text{\rm mod }\, p)
\\&=(D_{ij}(x^{(\al)}))^p{\otimes}1
+\sum\limits_{n,\ell=0}^{p{-}1} ({-}1)^{n+\ell} h(k,k')^{\lg
\ell\rg} h(m,m')^{\lg n \rg}{\otimes} (1{-}e(k,k')t)^{{-}\ell}\\ &
\quad \cdot (1{-}e(m,m')t)^{{-}n}
\Big(\delta_{\ell,0}\delta_{n,0}(D_{ij}(x^{(\al)}))^p-\delta_{n,0}\delta_{1,\ell}
(\delta_{ik}{-}\delta_{jk})\\
& \quad \cdot
\delta_{\alpha,\epsilon_i+\epsilon_j}e(k,k')-\delta_{\ell,0}\delta_{1,n}
(\delta_{im}{-}\delta_{jm})\delta_{\alpha,\epsilon_i+\epsilon_j}e(m,m')\Big
)t^{n+\ell}
 \\
&=(D_{ij}(x^{(\al)}))^p{\otimes}1+1{\otimes}(D_{ij}(x^{(\al)}))^p\\
& \quad
+h(k,k'){\otimes}(1{-}e(k,k')t)^{-1}(\delta_{ik}{-}\delta_{jk})
\delta_{\a,\epsilon_i+\epsilon_j}e(k,k')t\\
& \quad
+h(m,m'){\otimes}(1{-}e(m,m')t)^{-1}(\delta_{im}{-}\delta_{jm})
\delta_{\a,\epsilon_i+\epsilon_j}e(m,m')t.
\end{split}\tag{51}
\end{equation*}

Hence, when $\alpha\ne\epsilon_i+\epsilon_j$, we get
\begin{equation*}
\begin{split}
\Delta((D_{ij}(x^{(\al)}))^p)&\equiv(D_{ij}(x^{(\al)}))^p\otimes
1+1\otimes
(D_{ij}(x^{(\al)}))^p\\
&\in I_{t,q}\otimes
U_{t,q}(\mathbf{S}(n;\underline{1}))+U_{t,q}(\mathbf{S}(n;\underline{1}))\otimes
I_{t,q};
\end{split}
\end{equation*}
And when $\al=\epsilon_i+\epsilon_j$, by Lemmas 3.4 and 4.6, (45)
becomes
\begin{equation*}
\begin{split}
\Delta(D_{ij}(x^{(\epsilon_i+\epsilon_j)}))&=D_{ij}(x^{(\epsilon_i+\epsilon_j)})\otimes
1+ 1\otimes
D_{ij}(x^{(\epsilon_i+\epsilon_j)})\\
&\quad +\,(\delta_{ik}-\delta_{jk})h(k,k')\otimes
(1-e(k,k')t)^{-1}e(k,k')t\\
&\quad +\,(\delta_{im}-\delta_{jm})h(m,m')\otimes
(1-e(m,m')t)^{-1}e(m,m')t.
\end{split}
\end{equation*}
Combining with (51), we obtain
\begin{equation*}
\begin{split}
\Delta((D_{ij}(x^{(\epsilon_i+\epsilon_j)}))^p-D_{ij}(x^{(\epsilon_i+\epsilon_j)}))
&\equiv((D_{ij}(x^{(\epsilon_i+\epsilon_j)}))^p-D_{ij}(x^{(\epsilon_i+\epsilon_j)}))\otimes
1\\ &\quad +1\otimes ((D_{ij}(x^{(\al)}))^p-D_{ij}(x^{(\epsilon_i+\epsilon_j)}))\\
&\in I_{t,q}\otimes
U_{t,q}(\mathbf{S}(n;\underline{1}))+U_{t,q}(\mathbf{S}(n;\underline{1}))\otimes
I_{t,q}.
\end{split}
\end{equation*}

Thereby, we prove that the ideal $I_{t,q}$ is also a coideal of the
Hopf algebra $U_{t,q}(\mathbf{S}(n;\underline{1}))$.

\smallskip
(II) \ By Lemmas 4.5, 3.2, 3.4 \& 4.6, we have
\begin{equation*}
\begin{split}
S((D_{ij}(&x^{(\al)}))^p)
=-\bigl(1{-}e(k,k')t\bigr)^{-p\,\alpha(k,k')}\bigl(1{-}e(m,m')t\bigr)^{-p\,\alpha(m,m')}\\
& \quad \cdot \Bigl( \sum\limits_{n,\ell=0}^{\infty}
d_{m'}m^{(n)}d_{kk'}^{(\ell)}((D_{ij}(x^{(\al)}))^p)\cdot
h(k,k')_1^{\lg \ell\rg}h(m,m')_1^{\lg
n\rg}t^{n+\ell}\Bigr)\\
&\equiv -\sum\limits_{n,\ell=0}^{p-1}
d_{mm'}^{(n)}d_{kk'}^{(\ell)}((D_{ij}(x^{(\al)}))^p)\cdot
h(k,k')_1^{\lg \ell\rg}h(m,m')_1^{\lg n\rg}t^{n+\ell}
\quad (\text{mod }p)\\
&=-(D_{ij}(x^{(\al)}))^p+(\delta_{ik}{-}\delta_{jk})\delta_{\alpha,\epsilon_i+\epsilon_j}e(k,k')\cdot
h(k,k')_1^{\lg 1\rg} t\\
&\quad
+(\delta_{im}{-}\delta_{jm})\delta_{\alpha,\epsilon_i+\epsilon_j}e(m,m')\cdot
h(m,m')_1^{\lg 1\rg} t.
\end{split}\tag{52}
\end{equation*}

Hence, when $\alpha\ne\epsilon_i+\epsilon_j$, we get
$$
S\bigl((D_{ij}(x^{(\al)}))^p\bigr)=-(D_{ij}(x^{(\al)}))^p\in
I_{t,q};
$$
When $\al=\epsilon_i+\epsilon_j$, by Lemmas 3.4 and 4.5,
(46) reads as
$S(D_{ij}(x^{(\epsilon_i+\epsilon_j)}))=-D_{ij}(x^{(\epsilon_i+\epsilon_j)})+(\delta_{ik}-\delta_{jk})e(k,k')\cdot
h(k,k')_1^{\langle 1\rangle}t+(\delta_{im}-\delta_{jm})e(m,m')\cdot
h(m,m')_1^{\langle 1\rangle}t$. Combining with (52), we obtain
$$
S\bigl((D_{ij}(x^{(\epsilon_i+\epsilon_j)}))^p-D_{ij}(x^{(\epsilon_i+\epsilon_j)})\bigr)
=-\bigl((D_{ij}(x^{(\epsilon_i+\epsilon_j)}))^p-D_{ij}(x^{(\epsilon_i+\epsilon_j)})\bigr)
\in I_{t,q}.
$$

Thereby, we show that the ideal $I_{t,q}$ is indeed preserved by the
antipode $S$ of the quantization
$U_{t,q}(\mathbf{S}(n;\underline{1}))$.

\smallskip
(III) It is obvious to notice that
$\varepsilon((D_{ij}(x^{(\al)}))^p)=0$ for all $0\le\alpha\le\tau$.

\smallskip
This completes the proof.
\end{proof}
\begin{remark}
Corollary 4.2 gives more Drinfel'd twists. Using the same proof as
Theorem 4.7, we can get new families of noncommutative and
noncocommutative Hopf algebras of dimension $p^{1{+}(n-1)(p^n-1)}$
in characteristic $p$. Obviously, they are $p$-polynomial
deformations $(\mathbf{u}_{t,q}(\mathbf{S}(n;\underline{1})),m,
\iota,$ $\Delta,S,\varepsilon)$ of the restricted universal
enveloping algebra of $\mathbf{S}(n;\underline{1})$ over the
$p$-truncated polynomial ring $\mathcal{K}[t]_p^{(q)}$.
\end{remark}

\subsection{Different twisted structures} We shall show that the twisted structures given by
Drinfel'd twists with different product-length are nonisomorphic.
\begin{defi}
A Drinfel'd twist $\mathcal{F} \in A\otimes A$ on any Hopf algebra
$A$ is called \textit{compatible} if $\mathcal{F}$ commutes with the
coproduct $\Delta_0$.
\end{defi}
In other words, twisting a Hopf algebra $A$ with a
\textit{compatible} twist $\mathcal{F}$ gives exactly the same Hopf
structure, that is, $\Delta_{\mathcal{F}}=\Delta_0$. The set of
\textit{compatible} twists on $A$ thus forms a group.

\begin{lemm}$($\cite{MG}$)$
Let $\mathcal{F} \in A\otimes A$ be a  Drinfel'd twist on a Hopf
algebra $A$. Then the twisted structure induced by $\mathcal{F}$
coincides with the structure on $A$  if and only if $\mathcal{F}$ is
a compatible twist.
\end{lemm}

Using the same proof as in Theorem 4.1, we obtain

\begin{lemm}
Let $\mathcal{F}, \mathcal{G} \in A\otimes A$ be Drinfel'd twists on
a Hopf algebra $A$ with
$\mathcal{F}\mathcal{G}=\mathcal{G}\mathcal{F}$ and $\mathcal{F}\neq
\mathcal{G}$. Then $\mathcal{F}\mathcal{G}$ is a Drinfel'd twist.
 Furthermore, $\mathcal{G}$ is a Drinfel'd twist on
$A_{\mathcal{F}}$, $\mathcal{F}$ is a Drinfel'd twist on
$A_{\mathcal{G}}$
 and
$\Delta_{\mathcal{F}\mathcal{G}}=(\Delta_{\mathcal{F}})_{\mathcal{G}}
=(\Delta_{\mathcal{G}})_{\mathcal{F}}$.
\end{lemm}


Let $A$ denote one of objects:
$U(\mathbf{S}_{\mathbb{Z}}^{+})[[t]]$,
$U_{t,q}(\mathbf{S}(n,\underline{1}))$ and
$\mathbf{u}_{t,q}(\mathbf{S}(n,\underline{1}))$.
\begin{prop}
Drinfel'd twists
$\mathcal{F}^{\zeta(i)}:=\mathcal{F}(2,1)^{\zeta_1}\cdots\mathcal{F}(n,1)^{\zeta_{n-1}}$
$($where
$\zeta(i)=(\zeta_1,\cdots,\zeta_{n-1})=(\underbrace{1,\cdots,1}_i,0,{\cdots},0)\in
\mathbb Z_2^{n-1})$ lead to $n-1$ different twisted Hopf algebra
structures on $A$.
\end{prop}
\begin{proof} For $i=1$, $\mathcal{F}(2,1)$ gives one twisted structure with a
twisted coproduct different from the original one. For $i=2$, using
Lemma 4.11, we know that $\mathcal{F}(3,1)$ is a Drinfel'd twist and
not a compatible twist on
$U(\mathbf{S}_{\mathbb{Z}}^{+})[[t]]_{\mathcal{F}(2,1)}$. So the
twist  $\mathcal{F}(2,1)\mathcal{F}(3,1)$ gives new Hopf algebra
structure with the coproduct different from the previous one twisted
by $\mathcal{F}(2,1)$. Using the same discussion, we obtain that the
Drinfel'd twists $\mathcal{F}^{\zeta(i)}$ for
$\zeta(i)=(\underbrace{1,\cdots,1}_i,0,\cdots,0)\in
\mathbb{Z}_2^{n-1}$ give  $n{-}1$ different twisted structures on
$U(\mathbf{S}_{\mathbb{Z}}^{+})[[t]]$. This leads to the
corresponding result on $A$.
\end{proof}

\section{Quantizations of horizontal type for
$\mathbf{S}(n;\underline{1})$ and $\mathfrak{sl}_n$}

In this section, we assume that $n\geq 3$. Take $h:=\pa_k-\pa_{k'}$
and $e:=x^{\epsilon_k-\epsilon_m}\pa_m$ $(1\leq k \neq k'\neq m \leq
n)$ and denote by $\mathcal{F}(k,k';m)$ the corresponding Drinfel'd
twist. These twists will lead to the quantizations in horizontal
direction. So we call them the Drinfeld twists in {\it horizontal}
(while those twists used in Sections 3, 4 are in {\it vertical}).
Using the horizontal Drinfeld twists and the same discussion in
Sections 2, 3, we obtain some new quantizations of horizontal type
for the universal enveloping algebra of the special algebra
$\mathbf{S}(n;\underline{1})$.
The twisted structures given by the twists $\mathcal{F}(k,k';m)$ on
subalgebra $\mathbf{S}(n;\underline{1})_0$ are the same as those on
the special linear Lie algebra $\mathfrak{sl}_n$ over a field
$\mathcal{K}$ with $\text{char}(\mathcal{K})=p$ derived by the
Jordanian twists $\mathcal{F}=\mathrm{exp}(h\otimes \sigma)$,
$\sigma=\mathrm{ln}(1-e)$ for some two-dimensional carrier
subalgebra $B(2)=\text{Span}_{\mathcal K}\{h, e\}$ discussed in
\cite{KL}, \cite{KLS}, etc.

\subsection{Quantizations of horizontal type of $\mathbf u(\mathbf{S}(n;\underline{1}))$}
From Lemma 2.2 and Theorem 2.4, we have
\begin{lemm}
Fix two distinguished elements $h:=\pa_k{-}\pa_{k'}$,
$e:=x^{\epsilon_k-\epsilon_m}\pa_m$ $(1\leq k \neq k'\neq m \leq
n)$, the corresponding horizontal quantization of
$U(\mathbf{S}^+_{\mathbb{Z}})$ over
$U(\mathbf{S}^+_{\mathbb{Z}})[[t]]$ by Drinfel'd twist
$\mathcal{F}(k,k';m)$ with the product undeformed is given by
\begin{gather*}
\Delta(x^{\alpha}\pa)=x^{\alpha}\pa\otimes
(1{-}et)^{\alpha_k{-}\a_{k'}}{+}\sum\limits_{\ell=0}^{\infty}{({-}1)}^\ell
\,h^{\lg \ell\rg}\otimes(1{-}et)^{{-}\ell}\cdot\tag{53}  \\
\qquad\qquad \cdot\,
x^{\alpha{+}\ell(\epsilon_k-\epsilon_m)}(A_\ell\pa-B_\ell\pa_m)t^\ell,\\
S(x^{\alpha}\pa)={-}(1{-}et)^{-(\alpha_k{-}\a_{k'})}\cdot\Bigl(\sum\limits_{\ell=0}^{\infty}
\,x^{\alpha{+}\ell(\epsilon_k-\epsilon_m)}(A_\ell\pa-B_\ell\pa_m)\cdot
h_1^{\lg \ell\rg}t^\ell\Bigr),
\tag{54}\\
\varepsilon(x^{\alpha}\pa)=0,  \tag{55}
\end{gather*}
where $\alpha-\eta \in \mathbb{Z}^n_+,\, \eta=-\underline1,\,
A_\ell=\frac{1}{\ell!}\prod\limits_{j=0}^{\ell-1}(\al_m{-}j),\,
B_\ell=\pa(\epsilon_k-\epsilon_m) A_{\ell{-}1}$, with a convention
$A_0=1, A_{-1}=0$.
\end{lemm}

Note that $A_\ell=0$ for $\ell>\a_m$ and $B_\ell=0$ for $\ell>\a_m
+1$ in Lemma 5.1.
\begin{remark} According to the parametrization of the twists $\mathcal F(k,k';m)$, we
get $n(n{-}1)(n{-}2)$ {\it basic Drinfel'd twists} over
$U(\mathbf{S}^+_{\mathbb{Z}})$ and  consider the products of some
{\it basic Drinfel'd twists}. Using the same argument as Section 4,
one can get many more new Drinfel'd twists, which will lead to new
complicated quantizations not only over the
$U(\mathbf{S}_{\mathbb{Z}}^+)[[t]]$, but over the
$\mathbf{u}_{t,q}(\mathbf{S}(n;\underline{1}))$ as well.
\end{remark}

We firstly make {\it the modulo $p$ reduction} for the quantizations
of $U(\mathbf{S}^{+}_\mathbb{Z})$ in Lemma 5.1 to yield the
horizontal quantizations of $U(\mathbf{S}(n;\underline{1}))$ over
$U_t(\mathbf{S}(n;\underline{1}))$.
\begin{theorem}
Fix distinguished elements
$h=D_{kk'}(x^{(\epsilon_k+\epsilon_{k'})})$,
$e=D_{mk}(x^{(2\epsilon_k)})$ $(1\leq k\neq k'\neq m\leq n)$, the
corresponding horizontal quantization of
$U(\mathbf{S}(n;\underline{1}))$ over
$U_t(\mathbf{S}(n;\underline{1}))$ with the product undeformed is
given by
\begin{gather*}
\Delta(D_{ij}(x^{(\alpha)}))=D_{ij}(x^{(\alpha)})\otimes
(1{-}et)^{\a(k,k')} +\sum\limits_{\ell=0}^{p{-}1}{({-}1)}^\ell
 h^{\lg \ell\rg}\otimes(1{-}et)^{{-}\ell}\cdot\tag{56}\\
\qquad\qquad\quad \cdot\,\Bigl(\bar{A}_\ell
D_{ij}(x^{(\alpha{+}\ell(\epsilon_k-\epsilon_m))})+\bar{B}_\ell
(\delta_{ik}D_{jm}-\delta_{jk}D_{im})(x^{(\alpha{+}(\ell-1)(\epsilon_k-\epsilon_{m}))})\Bigr)t^\ell,
\end{gather*}
\begin{gather*}
S(D_{ij}(x^{(\alpha)}))={-}(1{-}et)^{-\a(k,k')}
\cdot\sum\limits_{\ell=0}^{p{-}1}\Bigl(\bar{A}_\ell
D_{ij}(x^{(\alpha{+}\ell(\epsilon_k-\epsilon_m))})\tag{57}\\
\qquad\qquad\qquad\qquad\quad +\,\bar{B}_\ell
(\delta_{ik}D_{jm}-\delta_{jk}D_{im})(x^{(\alpha{+}(\ell-1)(\epsilon_k-\epsilon_{m}))})\Bigr)
\cdot h_1^{\lg \ell\rg}t^\ell,
\\
\varepsilon(D_{ij}(x^{(\alpha)}))=0, \tag{58}
\end{gather*}
where $0\le \alpha \le \tau$, $\a(k,k')=
\alpha_k{-}\delta_{ik}{-}\delta_{jk}{-}\a_{k'}{+}\delta_{ik'}{+}\delta_{jk'}$,
$\bar A_\ell\equiv\binom{\alpha_k{+}\ell}{\ell}\,(\text{\rm mod}
\,p)$ for $0\leq \ell\leq\a_m$, $\bar B_\ell\equiv
\binom{\alpha_k{+}\ell-1}{\ell-1}\,(\text{\rm mod}\,p)$ for $1\leq
\ell\leq\a_m{+}1$, and otherwise, $\bar A_\ell=\bar B_\ell=0$.
\end{theorem}
\begin{proof}
Note that the elements
$\sum_{i,\alpha}\frac{1}{\alpha!}a_{i,\alpha}x^{\alpha}D_i$ in
$\mathbf{S}^+_{\mathcal{K}}$ for $0\le\alpha\le\tau$ will be
identified with $\sum_{i,\alpha}a_{i,\alpha}x^{(\alpha)}D_i$ in
$\mathbf{S}(n;\underline{1})$ and those in $J_{\underline{1}}$
(given in Section 3.1) with $0$. Hence, by Lemma 5.1, we get
\begin{equation*}
\begin{split}
\Delta(D_{ij}(x^{(\alpha)}))&
=\frac{1}{\a!}\Delta(x^{\alpha-\epsilon_i-\epsilon_j}(\a_j\pa_i{-}\a_i\pa_j))\\
&=D_{ij}(x^{(\alpha)})\otimes
(1{-}et)^{\a(k,k')}{+}\sum\limits_{\ell=0}^{p-1}{({-}1)}^\ell
\,h^{\lg \ell\rg}\otimes (1{-}et)^{{-}\ell}\cdot\\
& \quad  \cdot
\frac{1}{\a!}x^{\alpha-\epsilon_i-\epsilon_j{+}\ell(\epsilon_k-\epsilon_m)}
\bigl(A_\ell(\a_j\pa_i{-}\a_i\pa_j)-B_\ell\pa_m\bigr)t^\ell,
\end{split}
\end{equation*}
where
$A_\ell=\frac{1}{\ell!}\prod\limits_{j=0}^{\ell-1}(\al_m{-}\delta_{im}{-}\delta_{jm}{-}j),\,
B_\ell=(\a_j\pa_i{-}\a_i\pa_j)(\epsilon_k{-}\epsilon_m)
A_{\ell{-}1}$.

Write
\begin{gather*}
(*)=\frac{1}{\a!}x^{\alpha-\epsilon_i-\epsilon_j{+}\ell(\epsilon_k-\epsilon_m)}
\bigl(A_\ell(\a_j\pa_i{-}\a_i\pa_j)-B_\ell\pa_m\bigr), \\
(**)=\bar{A}_\ell
D_{ij}(x^{(\alpha{+}\ell(\epsilon_k-\epsilon_m))})+\bar{B}_\ell
(\delta_{ik}D_{jm}{-}\delta_{jk}D_{im})(x^{(\alpha{+}(\ell-1)(\epsilon_k-\epsilon_{m}))}).
\end{gather*}

We claim that $(*)=(**)$.

The proof will be given in the following steps:

$(\text{\rm i})$ \ For $\delta_{im}+\delta_{jm}=1$, we have
\begin{equation*}
(*)=\begin{cases}
\frac{(\a_k+\ell)!}{\a_k!}\frac{(\a_m-\ell)!}{\a_m!}(A_\ell{+}A_{\ell-1})
D_{ij}(x^{(\alpha{+}\ell(\epsilon_k-\epsilon_m))}), & \text{\it for
} \,
0\le \ell\leq\a_m,\\
0, & \text{\it for } \quad \ell>\a_m.
\end{cases}
\end{equation*}
A simple calculation shows that
$\frac{(\a_k+\ell)!}{\a_k!}\frac{(\a_m-\ell)!}{\a_m!}(A_\ell{+}A_{\ell-1})
\equiv\binom{\alpha_k{+}\ell}{\ell} \ (\text{\rm mod} \; p)$, for
$0\leq \ell\leq\a_m$. So, $(*)=(**)$.

$(\text{\rm ii})$ \ For $\delta_{im}+\delta_{jm}=0$, we consider
three subcases:

If $\delta_{ik}=1$, we have
\begin{equation*}
(*)=\begin{cases}
\frac{(\a_k+\ell)!}{\a_k!}\frac{(\a_m-\ell)!}{\a_m!}A_\ell
D_{kj}(x^{(\alpha{+}\ell(\epsilon_k-\epsilon_m))}) & \\
\ +\,\frac{(\a_k+\ell-1)!}{\a_k!}\frac{(\a_m-(\ell-1))!} {\a_m!}
A_{\ell{-}1} D_{jm}(x^{(\alpha{+}(\ell-1)(\epsilon_k-\epsilon_m))}),
& \text{\it for } \ 0\leq \ell\leq\a_m{+}1,\\
0, & \text{\it for } \quad \ell>\a_m{+}1.
\end{cases}
\end{equation*}
A simple calculation indicates that for $0\leq \ell\leq\a_m{+}1$,
\begin{equation*}
\begin{split}
\frac{(\a_k{+}\ell)!}{\a_k!}\frac{(\a_m{-}\ell)!}{\a_m!}A_\ell
&\equiv
\binom{\alpha_k{+}\ell}{\ell}=\bar A_\ell \ (\text{\rm mod} \; p),\\
\frac{(\a_k{+}\ell{-}1)!}{\a_k!}\frac{(\a_m{-}(\ell{-}1))!}
{\a_m!}A_{\ell{-}1}&\equiv\binom{\alpha_k{+}\ell{-}1}{\ell{-}1}=\bar
B_\ell \ (\text{\rm mod} \;p).
\end{split}
\end{equation*}
So, $(*)=(**)$.

If  $\delta_{jk}=1$, we have
\begin{equation*}
(*)=\begin{cases}
\frac{(\a_k+\ell)!}{\a_k!}\frac{(\a_m-\ell)!}{\a_m!}A_\ell
D_{ik}(x^{(\alpha{+}\ell(\epsilon_k-\epsilon_m))}) & \\
\ -\frac{(\a_k+\ell-1)!}{\a_k!}\frac{(\a_m-(\ell-1))!}
{\a_m!}A_{\ell-1}
D_{im}(x^{(\alpha{+}(\ell-1)(\epsilon_k-\epsilon_m))}), & \text{\it
for } \ 0\leq \ell\leq\a_m{+}1,\\
0, & \text{\it for } \quad \ell>\a_m{+}1.
\end{cases}
\end{equation*}
A simple computation shows that
\begin{equation*}
\begin{split}
\frac{(\a_k{+}\ell)!}{\a_k!}\frac{(\a_m{-}\ell)!}{\a_m!}A_\ell
&\equiv \binom{\alpha_k{+}\ell}{\ell}=\bar A_\ell\,(\text{\rm mod}
\;p),
\quad \text{\it for } 0\leq \ell\leq\a_m,\\
\frac{(\a_k{+}\ell{-}1)!}{\a_k!}\frac{(\a_m{-}(\ell{-}1))!}{\a_m!}
A_{\ell-1}&\equiv\binom{\alpha_k{+}\ell{-}1}{\ell{-}1}=\bar B_\ell
\,(\text{\rm mod} \;p), \quad \text{\it for } 0\leq
\ell\leq\a_m{+}1.
\end{split}
\end{equation*}
So, $(*)=(**)$.

If $\delta_{ik}=\delta_{jk}=0$, we have
$(*)=\frac{(\a_k+\ell)!}{\a_k!}\frac{(\a_m-\ell)!}{\a_m!}A_\ell
D_{ij}(x^{(\alpha{+}\ell(\epsilon_k-\epsilon_m))})$, and
\begin{equation*}
\begin{split}
\frac{(\a_k{+}\ell)!}{\a_k!}\frac{(\a_m{-}\ell)!}{\a_m!}A_\ell
&\equiv \binom{\alpha_k{+}\ell}{\ell}=\bar A_\ell \ (\text{\rm mod}
\;p),
\quad\text{\it for } \ 0\leq \ell\leq\a_m, \\
\bar B_\ell &\equiv 0 \ (\text{\rm mod} \;p), \quad\text{\it for } \
0\leq \ell\leq\a_m{+}1.
\end{split}
\end{equation*}
So, $(*)=(**)$.

\smallskip
Therefore, we verify the formula (56).

\smallskip
Applying a similar argument to the antipode, we can get the formula
(57).

\smallskip
This completes the proof.
\end{proof}

To describe $\mathbf{u}_{t,q}(\mathbf{S}(n;\underline{1}))$
explicitly, we still need an auxiliary Lemma.

\begin{lemm} Denote $e=D_{mk}(x^{(2\epsilon_k)})$, $d^{(\ell)}=\frac{1}{\ell!}(\text{\rm
ad}\,e)^\ell$. Then

\smallskip
$(\text{\rm i})$ \ $d^{(\ell)}(D_{ij}(x^{(\alpha)}))=\bar{A}_\ell
D_{ij}(x^{(\alpha{+}\ell(\epsilon_k-\epsilon_m))})$

\smallskip
\hskip3.2cm $+\,\bar{B}_\ell
(\delta_{ik}D_{jm}{-}\delta_{jk}D_{im})(x^{(\alpha{+}(\ell-1)(\epsilon_k-\epsilon_{m}))})$,

\smallskip
\hskip0.6cm where $\bar A_\ell, \bar B_\ell$ as in Theorem 5.3.

\smallskip
$(\text{\rm ii})$ \ $d^{(\ell)}(D_{ij}(x^{(\epsilon_i+\epsilon_j)}))
=\delta_{\ell,0}D_{ij}(x^{(\epsilon_i+\epsilon_j)})
-\delta_{1,\ell}(\delta_{ik}{-}\delta_{jk}{-}\delta_{im}{+}\delta_{jm})e$.

\smallskip
$(\text{\rm iii})$ \
$d^{(\ell)}((D_{ij}(x^{(\alpha)}))^p)=\delta_{\ell,0}(D_{ij}(x^{(\alpha)}))^p-\delta_{1,\ell}
(\delta_{ik}{-}\delta_{jk}{-}\delta_{im}{+}\delta_{jm})\delta_{\alpha,\epsilon_i+\epsilon_j}e$.
\end{lemm}

\begin{proof} We can get (i) from the proof of  Theorem 5.3.

(ii) Note that $\bar A_0=1$, $\bar B_0=0$. Using Theorem 5.3,
 for $\delta_{im}+\delta_{jm}=1$, we obtain $\bar A_1=1$ and  $\bar B_1=0$;
 for $\delta_{im}+\delta_{jm}=0$, we obtain $\bar A_1=0$ and  $\bar B_1=1$. We have
   $\bar A_\ell=\bar B_\ell=0$ for
$\ell>1$. Therefore, in any case, we arrive at the result as
desired.

\smallskip
(iii) From (15), we obtain that for $0\le \alpha\le\tau$,
\begin{equation*}
\begin{split}
d^{(1)}\,((D_{ij}(x^{(\alpha)}))^p)&=\frac{1}{(\alpha!)^p}\bigl[\,e,(D_{ij}(x^{\alpha}))^p\,\bigr]
=\frac{1}{(\alpha!)^p}\bigl[\,e,(x^{\alpha-\epsilon_i-\epsilon_j}(\al_j\partial_i{-}\al_i\pa_j))^p\,\bigr]\\
&=\frac{1}{(\alpha!)^p}\sum\limits_{\ell=1}^p(-1)^\ell\dbinom{p}
{\ell}(x^{\alpha-\epsilon_i-\epsilon_j}(\al_j\partial_i{-}\al_i\pa_j))^{p-\ell}\cdot \\
& \qquad\qquad \cdot
x^{\epsilon_k-\epsilon_m{+}\ell(\alpha{-}\epsilon_i-\epsilon_j)}
(a_\ell
\partial_m-b_\ell(\al_j\partial_i{-}\al_i\pa_j))\\
&\equiv
-\frac{a_p}{\alpha!}\,x^{\epsilon_k-\epsilon_m{+}p(\alpha{-}\epsilon_i-\epsilon_j)}
\pa_m\qquad(\text{mod }\,p\,)\\
&\equiv \begin{cases} -{a_p}\,e,\qquad & \text{\it if }\quad \alpha=\epsilon_i+\epsilon_j\\
0,\qquad & \text{\it if }\quad \alpha\ne\epsilon_i+\epsilon_j
\end{cases}\qquad(\text{mod }\,J),
\end{split}
\end{equation*}
where the last ``$\equiv$" by using the identification with respect
to modulo the ideal $J$ as before, and
$a_\ell=\prod\limits_{m=0}^{\ell-1}(\al_j\partial_i-\al_i\pa_j)(\epsilon_k-\epsilon_m+m(\al{-}\epsilon_i-\epsilon_j)),\
b_\ell=\ell\,\pa_m(\al{-}\epsilon_i-\epsilon_j)a_{\ell-1}$, and
$a_p=\delta_{ik}{-}\delta_{jk}{-}\delta_{im}{+}\delta_{jm}$, for
$\alpha=\epsilon_i+\epsilon_j$.

Consequently, by the definition of $d^{(\ell)}$, we get
$d^{(\ell)}((x^{(\alpha)}D_i)^p)=0$ in
$\mathbf{u}(\mathbf{S}(n;\underline{1}))$ for $2\le \ell\le p-1$ and
$0\le\alpha\le\tau$.
\end{proof}

Based on Theorem 5.3 and Lemma 5.4, we arrive at
\begin{theorem} Fix distinguished elements
$h=D_{kk'}(x^{(\epsilon_k+\epsilon_{k'})})$,
$e=D_{mk}(x^{(2\epsilon_k)})$ $(1\leq k\neq k'\neq m\leq n)$, there
exists a noncommutative and noncocummtative Hopf algebra $($of
horizontal type$)$
$(\mathbf{u}_{t,q}(\mathbf{S}(n;\underline{1})),m,\iota,\Delta,S,\varepsilon)$
over $\mathcal{K}[t]_p^{(q)}$ with the product undeformed, whose
coalgebra structure is given by
\begin{gather*}
\Delta(D_{ij}(x^{(\alpha)}))=D_{ij}(x^{(\alpha)})\otimes
(1{-}et)^{\a(k,k')} \tag{59}\\
\hskip5cm +\,\sum\limits_{\ell=0}^{p{-}1}{({-}1)}^\ell h^{\lg
\ell\rg}\otimes(1{-}et)^{{-}\ell}d^{(\ell)}(D_{ij}(x^{(\alpha)}))t^\ell,
\\
 S(D_{ij}(x^{(\alpha)}))={-}(1{-}et)^{-\a(k,k')}
\sum\limits_{\ell=0}^{p{-}1}d^{(\ell)}(D_{ij}(x^{(\alpha)}))
\cdot h_1^{\lg \ell\rg}t^\ell, \tag{60}\\
\varepsilon(D_{ij}(x^{(\alpha)}))=0, \tag{61}
\end{gather*}
where $0\le\alpha\le\tau$ and
$\a(k,k')=\alpha_k{-}\delta_{ik}{-}\delta_{jk}-\a_{k'}{+}\delta_{ik'}{+}\delta_{jk'}$,
which is finite dimensional with
$\dim_{\mathcal{K}}\mathbf{u}_{t,q}(\mathbf{S}(n;\underline{1}))=p^{1{+}(n-1)(p^n-1)}$.
\end{theorem}
\begin{proof} Utilizing the same arguments as in the proofs of Theorems 3.5 \& 4.7,
we shall show that the ideal
$I_{t,q}$ is a Hopf ideal of the {\it twisted} Hopf algebra
$U_{t,q}(\mathbf{S}(n;\underline{1}))$ as in Theorem 5.3. To this
end, it suffices to verify that $\Delta$ and $S$ preserve the
generators in $I_{t,q}$.

\smallskip
(I) \ By Lemmas 2.5, 5.3 \& 5.4 (iii), we obtain
\begin{equation*}
\begin{split}
\Delta(&(D_{ij}(x^{(\a)}))^p)=(D_{ij}(x^{(\a)}))^p\otimes
(1{-}et)^{p\,(\alpha_k{-}\a_{k'})}
\\
&\quad+\sum\limits_{\ell=0}^{\infty} ({-}1)^\ell h^{\lg
\ell\rg}\otimes
(1{-}et)^{{-}\ell}d^{(\ell)}((D_{ij}(x^{(\a)}))^p)t^\ell
\\
&\equiv(D_{ij}(x^{(\a)}))^p\otimes1+\sum\limits_{\ell=0}^{p{-}1}
({-}1)^\ell h^{\lg \ell\rg}\otimes
(1{-}et)^{{-}\ell}d^{(\ell)}((D_{ij}(x^{(\a)}))^p)t^\ell\quad
(\text{\rm mod }\, p)
\\
&=(D_{ij}(x^{(\a)}))^p\otimes1+1\otimes(D_{ij}(x^{(\a)}))^p \\
&\quad
+\,h\otimes(1{-}et)^{-1}(\delta_{ik}{-}\delta_{jk}{-}\delta_{im}
{+}\delta_{jm})\delta_{\a,\epsilon_i+\epsilon_j}et.
\end{split}\tag{62}
\end{equation*}
Hence, when $\alpha\ne\epsilon_i+\epsilon_j$, we get
\begin{equation*}
\begin{split}
\Delta((D_{ij}(x^{(\a)}))^p)&\equiv(D_{ij}(x^{(\a)}))^p\otimes
1+1\otimes
(D_{ij}(x^{(\a)}))^p\\
&\in I_{t,q}\otimes
U_{t,q}(\mathbf{S}(n;\underline{1}))+U_{t,q}(\mathbf{S}(n;\underline{1}))\otimes
I_{t,q}.
\end{split}
\end{equation*}
When $\al=\epsilon_i+\epsilon_j$, by Lemma 5.4 (ii), (56) becomes
\begin{equation*}
\begin{split}
\Delta(D_{ij}(x^{(\epsilon_i+\epsilon_j)}))&=D_{ij}(x^{(\epsilon_i+\epsilon_j)})\otimes
1+ 1\otimes D_{ij}(x^{(\epsilon_i+\epsilon_j)})\\
&\quad+\,h\otimes
(1{-}et)^{-1}(\delta_{ik}{-}\delta_{jk}{-}\delta_{im}{+}\delta_{jm})et.
\end{split}
\end{equation*}
Combining with (62), we obtain
\begin{equation*}
\begin{split}
\Delta((D_{ij}(x^{(\epsilon_i+\epsilon_j)}))^p-D_{ij}(x^{(\epsilon_i+\epsilon_j)})
)&\equiv\bigl((D_{ij}(x^{(\epsilon_i+\epsilon_j)}))^p-D_{ij}(x^{(\epsilon_i+\epsilon_j)})
\bigr)\otimes 1\\ &\quad +1\otimes
\bigl((D_{ij}(x^{(\epsilon_i+\epsilon_j)}))^p-D_{ij}(x^{(\epsilon_i+\epsilon_j)})
\bigr)\\
&\in I_{t,q}\otimes
U_{t,q}(\mathbf{S}(n;\underline{1}))+U_{t,q}(\mathbf{S}(n;\underline{1}))\otimes
I_{t,q}.
\end{split}
\end{equation*}

Thereby, we prove that $I_{t,q}$ is a coideal of the Hopf algebra
$U_{t,q}(\mathbf{S}(n;\underline{1}))$.

\smallskip
(II) \ By Lemmas 2.5, 5.3 \& 5.4 (iii), we have
\begin{equation*}
\begin{split}
S((D_{ij}(x^{(\a)}))^p) &=-(1{-}et)^{-p(\alpha_k-\a_{k'})}
\sum\limits_{\ell=0}^{\infty} d^{(\ell)}((D_{ij}(x^{(\a)}))^p)\cdot
h_1^{\lg
\ell\rg}t^\ell\\
&\equiv -(D_{ij}(x^{(\a)}))^p-\sum\limits_{\ell=1}^{p-1}
d^{(\ell)}((D_{ij}(x^{(\a)}))^p)\cdot h_1^{\lg \ell\rg}t^\ell \quad
(\text{mod } p)\\
&=-(D_{ij}(x^{(\a)}))^p
+(\delta_{ik}{-}\delta_{jk}{-}\delta_{im}{+}\delta_{jm})\delta_{\a,\epsilon_i+\epsilon_j}e\cdot
h_1^{\lg 1\rg} t.
\end{split}\tag{63}
\end{equation*}
Hence, when $\alpha\ne\epsilon_i+\epsilon_j$, we get
$$
S\bigl((D_{ij}(x^{(\a)}))^p\bigr)=-(D_{ij}(x^{(\a)}))^p\in I_{t,q}.
$$
When $\al=\epsilon_i+\epsilon_j$, by Lemma 5.4 (ii), (57) reads as
$$
S(D_{ij}(x^{(\epsilon_i+\epsilon_j)}))
=-D_{ij}(x^{(\epsilon_i+\epsilon_j)})+(\delta_{ik}{-}\delta_{jk}{-}
\delta_{im}{-}\delta_{jm}) e\cdot h_1^{\lg 1\rg} t.
$$
Combining with (63), we obtain
$$
S\bigl((D_{ij}(x^{(\epsilon_i+\epsilon_j)}))^p-D_{ij}(x^{(\epsilon_i+\epsilon_j)})\bigr)
=-\bigl((D_{ij}(x^{(\epsilon_i+\epsilon_j)}))^p-D_{ij}(x^{(\epsilon_i+\epsilon_j)})\bigr)
\in I_{t,q}.
$$

Thereby, we show that $I_{t,q}$ is preserved by the antipode $S$ of
$U_{t,q}(\mathbf{S}(n;\underline{1}))$ as in Theorem 5.3.

\smallskip
(III) It is obvious to notice that
$\varepsilon((D_{ij}(x^{(\a)}))^p)=0$ for all $0\le\alpha\le\tau$.

\smallskip
So, $I_{t,q}$ is a Hopf ideal in
 $U_{t,q}(\mathbf{S}(n;\underline{1}))$. We get a
finite-dimensional horizontal quantization on
 $\mathbf{u}_{t,q}(\mathbf{S}(n;\underline{1}))$.
\end{proof}

\subsection{Jordanian modular quantizations of
$\mathbf{u}(\mathfrak{sl}_n)$} Let $\mathbf{u}(\mathfrak{sl}_n)$
denote the restricted universal enveloping algebra of
$\mathfrak{sl}_n$. Since Drinfeld twists $\mathcal{F}(k,k';m)$ of
horizontal type closely act on the subalgebra $U((\mathbf S_{\mathbb
Z}^+)_0)[[t]]$, consequently on $\mathbf u_{t,q}(\mathbf
S(n;\underline 1)_0)$, these induce the Jordanian quantizations on
$\mathbf u_{t,q}(\mathfrak{sl}_n)$.

By Lemma 5.4 (i), we have
\begin{equation*}
\begin{split}
d^{(\ell)}(D_{ij}(x^{(2\epsilon_j)}))
&=\delta_{\ell,0}D_{ij}(x^{(2\epsilon_j)})+\delta_{1,\ell}
\bigl(\delta_{jm}D_{ik}(x^{(2\epsilon_k)})\\
&\quad -\delta_{ik}D_{mj}(x^{(2\epsilon_j)})+\delta_{jm}\delta_{ik}
D_{km}(x^{(\epsilon_k+\epsilon_m)})\bigr)-\delta_{2,\ell}\delta_{jm}\delta_{ik}e.
\end{split}
\end{equation*}
By Theorem 5.3, we have

\begin{theorem} Fix distinguished elements
$h=D_{kk'}(x^{(\epsilon_k+\epsilon_{k'})})$,
$e=D_{mk}(x^{(2\epsilon_k)})$ $(1\leq k\neq k'\neq m\leq n)$, the
corresponding Jordanian quantization of
$\mathbf{u}(\mathbf{S}(n,\underline{1})_0)\cong \mathbf
u(\mathfrak{sl}_n)$ over
$\mathbf{u}_{t,q}(\mathbf{S}(n,\underline{1})_0)\cong \mathbf
u_{t,q}(\mathfrak{sl}_n)$ with the product undeformed, whose
coalgebra structure is given by
\begin{gather*}
\Delta(D_{ij}(x^{(\epsilon_i+\epsilon_j)}))=D_{ij}(x^{(\epsilon_i+\epsilon_j)})\otimes
1+1\otimes D_{ij}(x^{(\epsilon_i+\epsilon_j)}) \tag{64}\\
\qquad\qquad\qquad\qquad
+\,(\delta_{ik}{-}\delta_{jk}{-}\delta_{im}{+}\delta_{jm})\,
h\otimes(1{-}et)^{-1}et,
\\
\Delta(D_{ij}(x^{(2\epsilon_j)}))=D_{ij}(x^{(2\epsilon_j)})\otimes
(1{-}et)^{\delta_{jk}{-}\delta_{ik}{-}\delta_{jk'}{+}\delta_{ik'}}
+1\otimes D_{ij}(x^{(2\epsilon_j)}) \tag{65}\\
\qquad\quad
-\,h\otimes(1{-}et)^{-1}\bigl(\delta_{jm}D_{ik}(x^{(2\epsilon_k)})-\delta_{ik}D_{mj}(x^{(2\epsilon_j)})+\delta_{jm}\delta_{ik}
D_{km}(x^{(\epsilon_k+\epsilon_m)})\bigr)t\\
\qquad -\,\delta_{jm}\delta_{ik}\,
h^{\lg 2\rg}\otimes(1{-}et)^{-2}et^2,\\
S(D_{ij}(x^{(\epsilon_i+\epsilon_j)}))=-D_{ij}(x^{(\epsilon_i+\epsilon_j)})+(\delta_{ik}
{-}\delta_{jk}{-}\delta_{im}{+}\delta_{jm})eh_1t, \tag{66} \\
S(D_{ij}(x^{(2\epsilon_j)}))=-(1{-}et)^{-(\delta_{jk}
{-}\delta_{ik}{-}\delta_{jk'}{+}\delta_{ik'})}
\cdot\Bigl(D_{ij}(x^{(2\epsilon_j)})+\tag{67} \\
\quad \bigl(\delta_{jm}D_{ik}(x^{(2\epsilon_k)})-\delta_{ik}D_{mj}
(x^{(2\epsilon_j)})+\delta_{jm}\delta_{ik}
D_{km}(x^{(\epsilon_k+\epsilon_m)})\bigr)h_1t-\delta_{jm}\delta_{ik}
eh_1^{\lg 2\rg}t^2\Bigr),  \\
\varepsilon(D_{ij}(x^{(\epsilon_i+\epsilon_j)}))=\varepsilon(D_{ij}(x^{(2\epsilon_j)}))=0.
\tag{68}
\end{gather*}
for $1\leq i\neq j\leq n$.
\end{theorem}

\begin{remark}
Since $\mathbf{S}(n,\underline{1})_0\cong \mathfrak{sl}_n$, which,
via the identification $D_{ij}(x^{(\epsilon_i+\epsilon_j)})$
 with $E_{ii}-E_{jj}$ and $D_{ij}(x^{(2\epsilon_j)})$ with $E_{ji}$
for $1\leq i\neq j\leq n$, we get a Jordanian quantization for
$\mathfrak{sl}_n$, which has been discussed by Kulish et al (cf.
\cite{KL}, \cite{KLS} etc.).
\end{remark}

\begin{coro} Fix distinguished elements
$h=E_{kk}-E_{k'k'}$, $e=E_{km}$ $(1\leq k\neq k'\neq m\leq n)$, the
corresponding Jordanian quantization of
$\mathbf{u}(\mathfrak{sl}_n)$ over
$\mathbf{u}_{t,q}(\mathfrak{sl}_n)$ with the product undeformed,
whose coalgebra structure is given by
\begin{gather*}
\Delta(E_{ii}-E_{jj})=(E_{ii}-E_{jj})\otimes 1+1\otimes(
E_{ii}-E_{jj})\tag{69}\\
\qquad\qquad\qquad \qquad +\,(\delta_{ik}{-}\delta_{jk}
{-}\delta_{im}{+}\delta_{jm})
h\otimes(1{-}et)^{-1}et,\\
\Delta(E_{ji})=E_{ji}\otimes
(1{-}et)^{\delta_{jk}{-}\delta_{ik}{-}\delta_{jk'}{+}\delta_{ik'}}
+1\otimes E_{ji}\tag{70}\\
\qquad\qquad\qquad\qquad
-\,h\otimes(1{-}et)^{-1}\bigl(\delta_{jm}E_{ki}
-\delta_{ik}E_{jm}\bigr)t-\delta_{jm}\delta_{ik}
h^{\lg 2\rg}\otimes(1{-}et)^{-2}et^2,\\
S( E_{ii}-E_{jj})=-(E_{ii}-E_{jj})+(\delta_{ik}
{-}\delta_{jk}{-}\delta_{im}{+}\delta_{jm})eh_1t, \tag{71}
\end{gather*}
\begin{gather*}
S(E_{ji})=-(1{-}et)^{-(\delta_{jk}{-}\delta_{ik}{-}\delta_{jk'}{+}\delta_{ik'})}
\Bigl(E_{ji}+\bigl(\delta_{jm}E_{ki}-\delta_{ik}E_{jm}\bigr)h_1t\tag{72}
\\
\qquad\qquad\qquad\qquad\qquad\qquad\qquad\qquad\quad
-\,\delta_{jm}\delta_{ik}
eh_1^{\lg 2\rg}t^2\Bigr), \\
\varepsilon(E_{ii}-E_{jj})=\varepsilon(E_{ji})=0. \tag{73}
\end{gather*}
for $1\leq i\neq j\leq n$.
\end{coro}
\begin{example} For $n=3$, take $h=E_{11}{-}E_{22}$, $e=E_{13}$, and set
$h'=E_{22}{-}E_{33}$, $f=(1{-}et)^{-1}$. By Corollary 5.8, we get a
Jordanian quantization on $\mathbf u_{t,q}(\mathfrak {sl}_3)$ with
the coproduct as follows (here we omit the antipode formulae which
can be directly written down from (71) \& (72)):
\begin{equation*}
\begin{split}
\Delta(h)&=h\otimes f+1\otimes h,\\
\Delta(h')&=h\otimes f+(h'{-}h)\otimes 1+1\otimes h',\\
\Delta(e)&=e\otimes f^{-1}+1\otimes e,\\
\Delta(E_{12})&=E_{12}\otimes f^{-2}+1\otimes E_{12},\\
\Delta(E_{21})&=E_{21}\otimes f^2+(1{+}h)\otimes E_{21}-h\otimes
fE_{21}f^{-1},\\
\Delta(E_{31})&=E_{31}\otimes f+(1{+}h)\otimes E_{31}-h\otimes fE_{31}f^{-1}+2(f^{-1}{-}1)E_{31}\otimes f(f{-}1),\\
\Delta(E_{23})&=E_{23}\otimes f+1\otimes E_{23},\\
\Delta(E_{32})&=E_{32}\otimes f^{-1}+(1{+}h)\otimes E_{32}-h\otimes
fE_{32}f^{-1},
\end{split}
\end{equation*}
where $\{f, h\}$ satisfying the relations: $[h,f]=f^2-f$, $h^p=h$,
$f^p=1$ generates the (finite-dimensional) Radford Hopf subalgebra
(with $f$ as a group-like element) over a field of characteristic
$p$.
\end{example}

\vskip10pt \centerline{\bf ACKNOWLEDGMENT}

\vskip10pt Authors are indebted to B. Enriquez and C. Kassel for
their valuable comments on quantizations when N.H. visited IRMA as
an invited professor of the ULP at Strasbourg from November to
December of 2007. N.H. is grateful to Kassel for his kind invitation
and extremely hospitality.

\bibliographystyle{amsalpha}

\begin{thebibliography}{A}
\bibitem {CP}V. Chari and A. Pressley, \textit{A Guide to Quantum Groups}, Cambridge
University Press, Cambridge, 1995.

\bibitem{DZ}D. Dokovic and K. Zhao, \textit{Derivations, isomorphisms and
second cohomology of generalized-Witt algebras}, Trans. Amer. Math.
Soc. \textbf{350} (1998), 643--664.

\bibitem{DZ1}D. Dokovic and K. Zhao, \textit{Generalized Cartan type S Lie
 algebras in characterisric zero}, J. Algebra, \textbf{193}
(1997), 144--179.

\bibitem {D}V.G. Drinfeld, \textit{Quantum groups}, Proceedings ICM
(Berkeley 1986) \textbf{1} (1987), AMS, 798--820.

\bibitem{EH} B. Enriquez and G. Halbout, \textit{Quantization of coboundary
Lie bialgebras}, arXiv Math: QA/0603740.

\bibitem{EK1} P. Etingof and D. Kazhdan, \textit{Quantization of
Lie bialgebras I}, Selecta Math. (N.S.), \textbf{2} (1) (1996),
1--41,  arXiv Math: q-alg/9506005.

\bibitem{EK2} P. Etingof and D. Kazhdan, \textit{Quantization of
Lie bialgebras II}, Selecta Math. (N.S.), \textbf{4} (2) (1998),
213--231, arXiv Math: q-alg/9701038.

\bibitem {ES} P. Etingof and O. Schiffmann, \textit{Lectures on Quantum Groups}, 2nd,
International Press, USA, 2002.

\bibitem{AJ}A. Giaquinto and J. Zhang, \textit{Bialgebra action, twists and
universal deformation formulas}, J. Pure Appl. Algebra, \textbf{128}
(2) (1998), 133--151.

\bibitem{MG}M. D. Gould and T. Lekatsas, \textit{Some twisted
results}, arXiv Math: QA/0504184.

\bibitem{CG}C. Grunspan, \textit{Quantizations of the Witt algebra and of simple
Lie algebras in characteristic $p$}, J. Algebra, \textbf{280}
(2004), 145--161.

\bibitem{HW}N. Hu and X. Wang, \textit{Quantizations of the generalized-Witt algebra and of Jacobson-Witt
 algebra in the modular case}, arXiv Math: QA/0602281, J. Algebra, \textbf{312} (2007),
 902--929.

\bibitem{KL}P. P. Kulish, V. D. Lyakhovsky and M. A. Olmo, \textit{Chains of twists for classical
Lie algebras}, arXiv Math: QA/9908061.

\bibitem{KLS}P. P. Kulish, V. D. Lyakhovsky and M. E. Samsonov,
\textit{Twists in $U(\mathfrak{sl}_3)$ and their quantizations},
arXiv Math: QA/0605392.

\bibitem{M} W. Michaelis, \textit{A class of infinite-dimensional Lie bialgebras containing
the Virasoro algebras}, Adv. Math., \textbf{107} (1994), 365--392.

\bibitem{NT} S. H. Ng and E.J. Taft, \textit{Classification of the Lie bialgebra structures
on the Witt and Virasoro algebras}, J. Pure and Appl. Algebra,
\textbf{151} (2000), 67--88.

\bibitem{VO}O. V. Ogievetsky, \textit{Hopf structures on the Borel subalgebra of ${\frak {sl}}(2)$},
Rend. Cir. Palermo (2) Suppl. \textbf{37} (1994), 185--199.

\bibitem{DR}D. E. Radford, \textit{Operators on Hopf algebras}, Amer. J. Math.
\textbf{99} (1977), 139--158.

\bibitem{NR}N. Reshetikhin, \textit{Multiparameter quantum groups and
twisted quasitriangular Hopf algebras}, Lett. Math. Phys.
\textbf{20} (1990), 331--335.

\bibitem{GY}G. Song and Y. Su, \textit{Lie
bialgeras of generalized-Witt type}, arXiv Math: QA/0504168, Science
in China, Ser. A---Math. \textbf{49} (4) (2006), 533--544.

\bibitem{R}R.P. Stanley, \textit{Enumerative Combinatorics, I},
Cambridge Studies in Advanced Mathematics, {\bf 49}, Cambridge
University Press, 1997.

\bibitem{H}H. Strade, \textit{Simple Lie Algebras over Fields of Positive Characteristic,
I. Structure Theory}, de Gruyter Expositions in Mathematics, {\bf
38}, Walter de Gruyter, 2004.

\bibitem{HR}H. Strade and R. Farnsteiner, \textit{Modular Lie Algebras and
Their Representations}, Monogr. Textbooks, Pure Appl. Math. {\bf
116}, Marcel Dekker, 1988.

\bibitem{T}E. Taft, \textit{Witt and Virasoro algebras as bialgebras}, J. Pure Appl.
Algebra \textbf{87} (3) (1993), 301--312.
\end{thebibliography}

\end{document}